\documentclass[review]{elsarticle}

\usepackage{lineno,hyperref}
\modulolinenumbers[5]

\journal{XXX}

\usepackage[colorinlistoftodos]{todonotes}
\usepackage[linesnumbered,ruled,vlined]{algorithm2e}
\usepackage[nottoc]{tocbibind}
\usepackage{booktabs,caption}
\usepackage{setspace}
\usepackage{graphicx}
\usepackage{svg}
\usepackage{adjustbox}
\usepackage{subcaption}
\usepackage{booktabs}
\graphicspath{ {./images/} }
\SetKwInput{kwInit}{Initialization}
\SetKwRepeat{Do}{do}{while}%

\usepackage{amssymb}% http://ctan.org/pkg/amssymb
\usepackage{pifont}% http://ctan.org/pkg/pifont

\usepackage[utf8]{inputenc}
\usepackage[english]{babel}
\usepackage[margin=25mm]{geometry}
\usepackage{amsmath,amssymb,mathtools}
\usepackage{amsthm}
\usepackage{bbm}
\usepackage{romannum}
\usepackage{multirow}
\usepackage[flushleft]{threeparttable}
\newtheorem{definition}{Definition}

\theoremstyle{definition}

\usepackage{ulem}

\usepackage{color}
\colorlet{Mycolor1}{green!20!orange!80!}
\iftrue  %commande pour voir le texte en couleur 
\newcommand{\comments}[1]{\footnote{\textcolor{blue}{\textit{#1}}}}
\newcommand{\DW}[1]{\textcolor{black}{#1}}
\newcommand{\AL}[1]{\textcolor{black}{#1}}
\else
\newcommand{\comments}[1]{}
\newcommand{}[1]{#1}
\newcommand{}[1]{#1}
\fi

\biboptions{authoryear}

\begin{document}
\pagenumbering{arabic}
\begin{frontmatter}

% \begin{titlepage}
% \begin{center}
% \vspace*{1cm}

% \textbf{ \large Optimization-Informed Neural Networks: a deep learning approach for solving \DW{constrained} nonlinear optimization problems}

% \vspace{1.5cm}

% % Author names and affiliations
% Dawen Wu$^{a}$ (dawen.wu@centralesupelec.fr), Abdel Lisser$^a$ (abdel.lisser@l2s.centralesupelec.fr) \\

% \hspace{10pt}

% \begin{flushleft}
% \small  

% \vspace{1cm}
% \textbf{Corresponding Author:} \\
% Dawen Wu \\
% Address: Université Paris-Saclay, CNRS, CentraleSupélec, Laboratoire des signaux et systèmes, 91190, Gif-sur-Yvette, France. \\
% Tel: (+33) 750798387 \\
% Email: dawen.wu@centralesupelec.fr

% \end{flushleft}        
% \end{center}
% \end{titlepage}

\title{Optimization-Informed Neural Networks}
% : a deep learning approach for solving \DW{constrained} nonlinear optimization problems}

\author[mymainaddress]{Dawen Wu \corref{mycorrespondingauthor}}
\author[mymainaddress]{Abdel Lisser}

\cortext[mycorrespondingauthor]{Corresponding author}
\ead{dawen.wu@centralesupelec.fr, abdel.lisser@l2s.centralesupelec.fr}

\address[mymainaddress]{Université Paris-Saclay, CNRS, CentraleSupélec, Laboratoire des signaux et systèmes, 91190, Gif-sur-Yvette, France}

\begin{abstract}
\DW{Solving constrained nonlinear optimization problems (CNLPs) is a longstanding computational problem that arises in various fields, e.g., economics, computer science, and engineering.
We propose optimization-informed neural networks (OINN), a deep learning approach to solve CNLPs. 
By neurodynamic optimization methods, a CNLP is first reformulated as an initial value problem (IVP) involving an ordinary differential equation (ODE) system.
A neural network model is then used as an approximate state solution for this IVP, and the endpoint of the approximate state solution is a prediction to the CNLP.
We propose a novel training algorithm that directs the model to hold the best prediction during training. 
In a nutshell, OINN transforms a CNLP into a neural network training problem.
By doing so, we can solve CNLPs based on deep learning infrastructure only, without using standard optimization solvers or numerical integration solvers.
The effectiveness of the proposed approach is demonstrated through a collection of classical problems, e.g., variational inequalities, nonlinear complementary problems, and standard CNLPs. }
% Some challenging CNLPs seem difficult to solve by numerical methods.
% We compare the solution process between the OINN and numerical methods.
\end{abstract}

\begin{keyword}
Constrained nonlinear optimization problems \sep Neural networks \sep Neurodynamic optimization \sep ODE system 
\end{keyword}

\end{frontmatter}
% \linenumbers

\section{Introduction}
% NLPs description
Constrained nonlinear optimization problems (CNLPs) play a central role in operations research and \DW{have a wide range of real-world applications, such as production planning, resource allocation, portfolio selection, portfolio optimization, feature selection, equilibrium problems \citep{xiao2006optimal, leung2020minimax, wang2021neurodynamic, wu2022dynamical}.} CNLPs have been studied at both the theoretical and practical levels for the last few decades \citep{bertsekas1997nonlinear, boyd2004convex}.
% NLPs are a broad concept encompassing many different forms, such as nonlinear projection problems, variational inequalities, and standard NLP.

% Neurodynamic approaches
\DW{Neurodynamic optimization methods model a CNLP by the mean of an ordinary differential equation (ODE) system.  \cite{hopfield1985neural} pioneered this study and solved the well-known “traveling salesman” problem by the Hopfield network.} \cite{kennedy1988neural} extended the method to solve nonlinear convex programming problems by using a penalty parameter. 
However, the disadvantage of this penalty parameter method is that the true minimizer is obtained only when the penalty parameter goes to infinity. 
When the penalty parameter is too large, the method hardly converges to the optimal solution. 
Since then, researchers have improved the method gradually without using the penalty parameter. \cite{rodriguez1990nonlinear, xia2002projection, gao2004neural, xia2007new, xia2015bi} proposed neurodynamic methods based on a projection function. Besides the convex and smooth optimization problems, \cite{forti2004generalized, xue2008subgradient, qin2014two} solved non-smooth CNLPs using differential inclusion theory and subgradient. 
Additionally, pseudoconvex optimization problems have been studied based on various assumptions \citep{guo2011one, qin2013new, xu2020neurodynamic}. 

% NNs
With the rapid growth of available data and computing resources, deep learning now has a wide range of applications, e.g., image processing \citep{krizhevsky2012imagenet, Goodfellow-et-al-2016}, natural language processing \citep{DBLP:journals/corr/abs-1810-04805}, bioinformatics \citep{min2017deep, jumper2021highly}. In operations research, a neural network is used as a solver component to solve the mixed integer programming problem \citep{nair2020solving}.
Graphical neural networks can be used for combinatorial optimization problems directly as solvers or to enhance standard solvers \citep{cappart2021combinatorial}.

% ODEs-NN
\cite{https://doi.org/10.1002/cnm.1640100303} initially used a neural network as an approximate solution to differential equations, where the training objective is to satisfy the given differential equation and boundary conditions. 
\cite{712178} constructed a neural network to satisfy an initial/boundary condition, and they discussed the use of ODE and PDE problems, respectively. 
\cite{lagaris2000neural, mcfall2009artificial} extended the Lagris' method  to irregular boundaries.
\cite{RAISSI2019686} introduced physics-informed neural networks to solve forward and inverse problems involving PDEs. 
\DW{\cite{sirignano2018dgm} presents a theoretical analysis that shows the neural network approximator converges to the PDE solution as hidden units go to infinity.
Deep learning approaches are attempting to overcome the challenge of solving high-dimensional nonlinear PDEs \citep{han2017deep, yu2018deep, han2018solving, beck2019machine}.
This line of research has been extended to various fields, e.g.,  computational mechanics \citep{anitescu2019artificial, samaniego2020energy, guo2021deep}.}
All the above methods use one neural network to solve one ODE/PDE problem. \cite{flamant2020solving} parameterizes ODE systems and uses the parameters as an input to a neural network so that one neural network can solve multiple ODE systems.
The universal approximation theorem of neural networks states that a neural network can approximate any continuous function to arbitrary accuracy \citep{Cybenko1989, HORNIK1989359, SONODA2017233}. 
Automatic differentiation tools facilitate the computation of derivative, gradient, and Jacobian matrix \citep{baydin2018automatic, NEURIPS2019_9015}. Software packages have been developed to implement these deep learning methods for solving differential equations \citep{2021, Chen2020}.

\subsection{Contributions}
The contributions of this paper can be summarized as follows.
\begin{itemize}
    \item \DW{We propose a deep learning approach to solve CNLPs, namely OINN. To the best of our knowledge, this is the first time deep learning is used to solve CNLPs. OINN reformulates a CNLP as a neural network training problem via neurodynamic optimization. \AL{Thus,} we can solve the CNLP by only deep learning infrastructure without using any standard CNLP solvers or numerical integration solvers.}
    \item \DW{We propose a dedicated algorithm to train the OINN model toward solving the CNLP. This algorithm is based on the epsilon metric, which is used to evaluate approximate solutions to the CNLP.}
    \item 
    \DW{We present the difference between OINN and numerical integration methods for solving a CNLP.
    OINN can give an approximate solution at any round of iterations, while the numerical integration methods can only give the solution at the end of the program.
    We show the computational advantages of OINN thanks to this feature.}
\end{itemize}

\DW{The remaining sections are organized as follows. 
The background information necessary to understand this paper is provided in Section \ref{section: 2}, including an introduction of CNLPs, neurodynamic optimization methods, and numerical integration methods.
Section \ref{section: 3} describes the OINN model and how it solves a CNLP.
Training of the OINN model is given in Section \ref{section: 4}.
Experimental results \AL{are given in Section \ref{section: 5}}  where six different CNLP instances are solved with OINN, and a comparison between OINN and numerical integration methods is presented.
Section \ref{section: 6} summarizes this paper and gives future directions.}

\subsection{Notations}
The notations list for this paper is shown in Table \ref{tab: Notations}.
\begin{table}[h]
\begin{tabular}{@{}ll@{}}
\toprule
Notation                       & Definition                                                                                                                  \\ \midrule
$\mathbf{y} \in \mathbb{R}^n$                 & Variables of a CNLP. $n$ refers to the number of variables.\\
$\mathbf{y}^* \in \mathbb{R}^n$                 & Optimal solution of a CNLP.\\
$\mathbf{x} \in \mathbb{R}^{j}, \mathbf{u} \in \mathbb{R}^{k}$    & \begin{tabular}[l]{@{}l@{}}Primal variables and dual variables of a standard CNLP.\\ $j$ refers to the number of primary variables, $k$ refers to the number of dual variables\end{tabular}\\
$P(\cdot)$              & A projection function that map variables onto a feasible set.               \\
$\Phi(\mathbf{y}) = \frac{d\mathbf{y}}{dt}: \mathbb{R}^{n} \to \mathbb{R}^{n}$         & An ODE system            \\
\DW{$\bar{\mathbf{y}}(t): \mathbb{R} \to \mathbb{R}^n$}              & \DW{True state solution of an ODE system}                                                             \\
\DW{$\hat{\mathbf{y}}(t): \mathbb{R} \to \mathbb{R}^n$}              & \DW{Approximate state solution obtained by numerical integration methods}                                                             \\
$\mathbf{y}(t; \mathbf{w})$  & An OINN model, where $\mathbf{w}$ are the model parameters                                          \\
$\mathbf{y}(T; \mathbf{w}) \in \mathbb{R}^n$  & The endpoint of an OINN model                                          \\
$\mathbf{y}_0\in \mathbb{R}^n$        &  An initial point of an ODE system                                        \\
$[0, T] \subset \mathbb{R}$                   & A time range of an ODE system                                                     \\
$\epsilon(\cdot)$  & Epsilon metric for evaluating a solution of CNLP                                           \\
$\mathcal{L}(t; \mathbf{w})$ & Loss function of OINN            \\
$E(\mathbf{w})$ & Objective function of OINN            \\
OINN  & Optimization-informed neural networks                                                                 \\
CNLP  & Constrained nonlinear optimization problem                                                                          \\
ODE   & Ordinary differential equation                                                                 \\
IVP  & Initial value problem                                                                 \\
\DW{NPE}  & \DW{Nonlinear projection equation}                                
\\ \bottomrule
\end{tabular}
\centering
\caption{Notations}
\label{tab: Notations}
\end{table}

\section{Preliminaries}\label{section: 2}
Serval types of CNLP are introduced in Section \ref{section: 2.1}.
Neurodynamic optimization methods, which model a CNLP as an ODE system, are introduced in Section \ref{section: 2.2}.
\DW{The initial value problem and numerical integration methods are described in Section \ref{section: 2.3}.}

\subsection{Constrained nonlinear optimization problems}\label{section: 2.1}

This subsection introduces four types of CNLP, i.e., standard CNLP, variational inequality, nonlinear complementary problem, and nonlinear projection equation.

% Standard CNLP
\textbf{Standard CNLP} \ The standard CNLP has the following form
\begin{equation}
\begin{aligned}
\begin{cases}
&\min\limits_\mathbf{x} f(\mathbf{x}) \\
&\text{s.t.}           \\
&\quad  g(\mathbf{x})\leq \mathbf{0},    \\
&\quad  \mathbf{A} \mathbf{x} = \mathbf{b}, 
\end{cases}
\end{aligned}
\label{eq: CNLP}
\end{equation}
where $\mathbf{x} \in \mathbb{R}^j$ is the primal variable, $\mathbf{u} \in \mathbb{R}^k$ is the dual variable associated with the constraint $g(\mathbf{x})$. The objective function $f(\mathbf{x}): \mathbb{R}^{j}\to \mathbb{R}$ is not necessary convex or smooth. The constraint function $g(\mathbf{x}): \mathbb{R}^{j}\to \mathbb{R}^{k}$ is convex but not necessarily smooth, $\mathbf{A}\in \mathbb{R}^{e \times j}$ and  $\mathbf{b}\in \mathbb{R}^{e}$. 
The following projection function can project the variable $\mathbf{x}$ onto the equality constraints feasible set $\{\mathbf{x} \in \mathbb{R}^{j} \mid \mathbf{A} \mathbf{x}=\mathbf{b}\}$
\begin{equation}
P_{eq}(\mathbf{x})=\mathbf{x}-\mathbf{A}^{T}\left(\mathbf{A} \mathbf{A}^{T}\right)^{-1}(\mathbf{A} \mathbf{x}-\mathbf{b}).
\label{eq: linear constraint projection}
\end{equation}
The standard CNLP \eqref{eq: CNLP} is the most common form of CNLP, and we can classify it according to the property of the objective and constraint functions.
For example, it is called quadratic programming if the objective function is quadratic and the constraints are linear;
Nonsmooth optimization problems are those that involve non-smooth functions.
% The standard CNLP \eqref{eq: CNLP} can cover optimization models in portfolio optimization, optimal transportation,  financial engineering, power systems, etc \citep{bertsekas1997nonlinear, boyd2004convex}.

% Nonlinear projection equation
\textbf{Nonlinear projection equation (NPE)} \ A NPE aims at finding a vector $\mathbf{y}^*\in\mathbb{R}^n$ such that satisfies
\begin{equation}
P_{\Omega}(\mathbf{y}-G(\mathbf{y}))=\mathbf{y},
\label{eq: NPE}
\end{equation}
where $\mathbf{y} \in \mathbb{R}^n$ is a real vector, $G(\cdot): \mathbb{R}^n \to \mathbb{R}^n$ is a locally Lipschitz continuous function, $\Omega=\{\mathbf{y} \in \mathbb{R}^n | l_i^-\leq y_i\leq l_i^+, i=1,\dots,n\}$ is a box-constrained feasible set, where $l_i^-$ and $l_i^+$ are the lower and upper bounds of $y_i$, respectively. $P_\Omega(\cdot): \mathbb{R}^{n} \to \Omega$ is a projection function that project the variable onto the feasible set $\Omega$, defined by
\begin{equation}
P_{\Omega}(\mathbf{s})=\left(P^1_{\Omega}\left(s_{1}\right), \ldots, P^n_{\Omega}\left(s_{n}\right)\right)^{T}, 
\quad \text{where} \ 
P^i_{\Omega}\left(s_{i}\right)= \begin{cases}l_i^-, & \text{if}\  s_{i}<l_i^- \\ s_{i}, & \text{if} \ l_i^- \leq s_{i} \leq l_i^+ \\ l_i^+, & \text{otherwise}.\end{cases}
\label{eq: projection 1}
\end{equation}

% variational inequality
\textbf{Variational inequality (VI)} A VI aims at finding a vector $\mathbf{y}^*\in \Omega$ such that the following inequalities hold
\begin{equation}
\left(\mathbf{y}-\mathbf{y}^{*}\right)^{T} G\left(\mathbf{y}^{*}\right) \geq 0, \quad \mathbf{y} \in \Omega.
\label{eq: VI}
\end{equation}
where $G(\cdot)$ and $\Omega$ are the same as in \eqref{eq: NPE}. Variational inequality provides a reformulation of the Nash equilibrium in game theory to study equilibrium properties, including existence, uniqueness, and convergence \citep{patriksson2013nonlinear, parise2019variational, singh2018variational}.

% nonlinear complementary problem
\textbf{Nonlinear complementary problem (NCP)} A NCP is to find out a vector $\mathbf{y}^*$ that satisfies
\begin{equation}
G(\mathbf{y}) \geq 0, \quad \mathbf{y} \geq 0, \quad G(\mathbf{y})^{T} \mathbf{y}=0.
\label{eq: NCP}
\end{equation}
where $G(\cdot)$ is the same as in \eqref{eq: NPE}. Nonlinear complementarity problems arise in many practical applications. For example, finding a Nash equilibrium is a special case of the Linear complementarity problem; KKT systems of mathematical programming problems can be formulated as NCP problems.

Both the variational inequality \eqref{eq: VI} and nonlinear complementary problem \eqref{eq: NCP} can be reformulated as an NPE problem \citep{harker1990finite, robinson1992normal}. In addition, when the objective and constraint functions are convex and smooth, the CNLP \eqref{eq: CNLP} can also be reformulated as an NPE problem, where the variable vector is composed by the primal and dual variable, i.e., $\mathbf{y}=(\mathbf{x}^T, \mathbf{u}^T)^T$.

\subsection{Neurodynamic optimization}\label{section: 2.2}

This subsection introduces neurodynamic optimization methods, which model a CNLP by an ODE system. 
Consider a CNLP with an optimal solution $\mathbf{y}^{*}$. 
A neurodynamic approach establishes a dynamical system in the form of a first-order ODE system, i.e., $\frac{d\mathbf{y}}{dt} = \Phi(\mathbf{y})$. 
The state solution $\mathbf{y}(t)$ is expected to converge to the optimal solution of the CNLP, i.e., $\lim _{t \rightarrow \infty} \mathbf{y}(t)=\mathbf{y}^{*}$.  
Here, we present three different neurodynamic approaches \citep{xia2007new, qin2014two, xu2020neurodynamic}, each of which solves a type of CNLP.

\begin{definition}
Consider an ODE system $\frac{d\mathbf{y}}{dt}=\Phi(\mathbf{y})$, where $\Phi(\mathbf{y}): \mathbb{R}^n \to \mathbb{R}^n$. Given  a point $(t_0, \mathbf{y}_0)\in\mathbb{R}^{n+1}$, a vector value function $\mathbf{y}(t): \mathbb{R} \to \mathbb{R}^n$ is called a state solution, if it satisfies the ODE system $\frac{d\mathbf{y}}{dt}=\Phi(\mathbf{y})$ and the initial condition $\mathbf{y}(t_0) = \mathbf{y}_0$.
\label{definition: state solution}
\end{definition}

\cite{xia2007new} proposed a neurodynamic approach to model the nonlinear projection equation \eqref{eq: NPE}. The ODE system is as follows
\begin{equation}
    \frac{\mathrm{d} \mathbf{y}}{\mathrm{~d} t}=\lambda\left( - G\left(P_{\Omega}(\mathbf{y})\right)+P_{\Omega}(\mathbf{y})-\mathbf{y} \right),
\label{eq: 2007 XIA}
\end{equation}
where $\lambda > 0 $ is a parameter controlling the convergence rate.

\cite{qin2014two} proposed a neurodynamic approach to model the standard CNLP \eqref{eq: CNLP} when the objective function $f(\cdot)$ and the constraint function $g(\cdot)$ are convex and nonsmooth. The ODE system is as follows
\begin{equation}
    \begin{aligned} 
        \frac{d\mathbf{x}}{dt} \in &-(\mathbf{I}-\mathbf{U})\left[\partial f(\mathbf{x})+\partial g(\mathbf{x})^{T} (\mathbf{u}+g(\mathbf{x}))^{+} \right] - \mathbf{A}^{T} \rho(\mathbf{A} \mathbf{x}-\mathbf{b}), \\ 
        \frac{d\mathbf{u}}{dt} =& \frac{1}{2}\left(-\mathbf{u}+(\mathbf{u}+g(\mathbf{x}))^{+}\right),
    \end{aligned}
\label{eq: 2014 QIN}
\end{equation}
where $\mathbf{U}=\mathbf{A}^{T}\left(\mathbf{A} \mathbf{A}^{T}\right)^{-1} \mathbf{A}$, $\mathbf{I}$ is the identity matrix, and $\rho(\cdot)$ is defined as
\begin{equation}
\rho(\mathbf{s})=\left(\tilde{\rho}\left(s_{1}\right), \tilde{\rho}\left(s_{2}\right), \ldots, \tilde{\rho}\left(s_{e}\right)\right)^{T}, \quad \text{where} \ 
\tilde{\rho}\left(s_{i}\right)= \begin{cases}1, & \text { if } s_{i}>0 \\ {[-1,1],} & \text { if } s_{i}=0 \\ -1. & \text { otherwise } \end{cases}
\label{eq: h}
\end{equation}

\cite{xu2020neurodynamic} proposed a neurodynamic approach to model the standard CNLP \eqref{eq: CNLP} when the objective function $f(\cdot)$ is a pseudoconvex nonsmooth function, and the constraint function $g(\cdot)$ is a convex nonsmooth function. Unlike \eqref{eq: 2014 QIN}, the ODE system only models the state of $\mathbf{x}$ without considering $\mathbf{u}$. It is defined as follows
\begin{equation}
    \frac{d\mathbf{x}}{d t} \in-\theta(t)(\mathbf{I}-\mathbf{U})\left(\left\{\prod_{i=1}^{k}\left(1-\mu\left(g_{i}(\mathbf{x})\right)\right)\right\} \partial f(\mathbf{x}) + \partial B(\mathbf{x}) \right) - \mathbf{A}^{T} \rho(\mathbf{A} \mathbf{x}-\mathbf{b}),
\label{eq: 2020 QIN}
\end{equation}
where $U=\mathbf{A}^{T}\left(\mathbf{A} \mathbf{A}^{T}\right)^{-1} \mathbf{A}$, and $\rho(\cdot)$ is the same as in \eqref{eq: h}. $\theta(t)$ is defined by 
\begin{equation}
\theta(t)= \begin{cases}0, & \text { if } t \leq T_{0} \\ 1, & \text { otherwise } \end{cases}
\end{equation}
where $T_{0}=1+\left\|\mathbf{A} \mathbf{x}_{0}-b\right\|_{1} / \lambda_{\min }\left(\mathbf{A} \mathbf{A}^{\mathrm{T}}\right)$, $\mathbf{x}_0$ is an initial point, $\lambda_{\min }\left(\mathbf{A} \mathbf{A}^{\mathrm{T}}\right)=\min \left\{\lambda: \lambda\right.$ is the eigenvalue of $\left.\mathbf{A} \mathbf{A}^{\mathrm{T}}\right\}$. $\mu(\cdot)$ is defined by
\begin{equation}
\mu(s)= \begin{cases}1, & \text { if } s>0 \\ {[0,1],} & \text { if } s=0 \\ 0, & \text { if } s<0\end{cases}
\end{equation}
$\partial B(\mathbf{x})$ is given by
\begin{equation}
\partial B(\mathbf{x})= \begin{cases}\{0\}, & \mathbf{x} \in S \cap \operatorname{int}(\mathcal{F}) \\ \sum_{i \in I^{0}(\mathbf{x})} \mu\left(g_{i}(\mathbf{x})\right) \partial g_{i}(\mathbf{x}), & \mathbf{x} \in S \cap \operatorname{bd}(\mathcal{F}) \\ \sum_{i \in I^{0}(\mathbf{x})} \mu\left(g_{i}(\mathbf{x})\right) \partial g_{i}(\mathbf{x})+\sum_{i \in I^{+}(\mathbf{x})} \partial g_{i}(\mathbf{x}), & \mathbf{x} \in S \backslash \mathcal{F}\end{cases}
\end{equation}
where $\mathcal{F}=\left\{\mathbf{x}: g_{i}(\mathbf{x}) \leq 0, i=1,2, \ldots, k\right\}$, $S=\{x: \mathbf{A} x=\mathbf{b}\}$, $I^{0}(\mathbf{x})=\left\{i \in\{1,2, \ldots, k\}: g_{i}(\mathbf{x})=0\right\}$, $I^{+}(\mathbf{x})=\left\{i \in\{1,2, \ldots, k\}: g_{i}(\mathbf{x})>0\right\}$.

\begin{table}[t]
\adjustbox{width=\textwidth}{
\begin{tabular}{llllll}
\hline
\multirow{2}{*}{Reference} & \multicolumn{3}{l}{Neurodynamic approach}                                  & \multirow{2}{*}{CNLP}                                   & \multirow{2}{*}{Conditions on functions}                                                                                     \\ \cline{2-4}
                           & ODE system       & Projection                          & Global convergence &                                                                                                            &                                                                                                                              \\ \hline
\cite{xia2007new}                   & \eqref{eq: 2007 XIA} & \eqref{eq: projection 1}              & True       & \begin{tabular}[c]{@{}l@{}}\eqref{eq: CNLP}, \eqref{eq: NPE}\\ \eqref{eq: VI}, \eqref{eq: NCP}\end{tabular} & \begin{tabular}[c]{@{}l@{}} In \eqref{eq: CNLP}, $f(\cdot)$ and $g(\cdot)$ are \\ convex and smooth \\ In \eqref{eq: NPE}, \eqref{eq: VI}, \eqref{eq: NCP}, $G(\cdot)$ is locally lipschitz.\end{tabular} \\\hline
\cite{qin2014two}                  & \eqref{eq: 2014 QIN} & \eqref{eq: linear constraint projection} & True       & \eqref{eq: CNLP}                                                                                            & \begin{tabular}[c]{@{}l@{}}In \eqref{eq: CNLP}, $f(\cdot)$ and $g(\cdot)$ are \\ convex and nonsmooth\end{tabular}                                \\\hline
\cite{xu2020neurodynamic}                   & \eqref{eq: 2020 QIN} & \eqref{eq: linear constraint projection} & True      & \eqref{eq: CNLP}                                                                                            & \begin{tabular}[c]{@{}l@{}}In \eqref{eq: CNLP}, $f(\cdot)$ and $g(\cdot)$ are \\ pseudoconvex and smooth\end{tabular}                             \\ \hline
\end{tabular}}
\caption{\textbf{Summary of the neurodynamic approaches and their corresponding CNLPs.} Global convergence refers to whether the ODE system can globally converge to the solution set of the CNLP. Conditions on functions refer to the assumptions required for the CNLP.}
\centering
\label{tab: 3 problems}
\end{table}

\begin{definition}
An ODE system $\frac{d\mathbf{y}}{dt} = \Phi(\mathbf{y})$ is said to be globally converges to a solution set $\mathcal{Y}^*$ if for any given initial point, the state solution $\mathbf{y}(t)$ satisfies
$$\lim _{t \rightarrow \infty} \operatorname{dist}\left(\mathbf{y}(t), \mathcal{Y}^{*}\right)=0,$$
where $\operatorname{dist}\left(\mathbf{y}(t), \mathcal{Y}^{*}\right)=\inf _{\mathbf{y}^* \in \mathcal{Y}^{*}}\|\mathbf{y}(t)-\mathbf{y}^*\|$, and $\left\lVert \cdot \right\rVert$ is the euclidean norm. In particular, if the set $\mathcal{Y}^*$ contains only one point $\mathbf{y}^*$, then $\lim _{t \rightarrow \infty} \mathbf{y}(t)=\mathbf{y}^*$, and the ODE system is globally asymptotically stable at $\mathbf{y}^*$.
\label{definition: globally convergence}
\end{definition}

% Global convergence
The global convergence property states that starting from any initial point, the state solution $\mathbf{y}(t)$ of the ODE system converges to the CNLP solution as time $t$ goes to infinity. 
A neurodynamic approach usually establishes the global convergence property in two steps: First, the ODE system's equilibrium points coincide with the optimal solutions of the CNLP. Then, using Lyapunov's theorem or LaSalle's invariance principle to prove that any state solution will converge to an equilibrium point of the ODE system.

% Description of the table
Table \ref{tab: 3 problems} summarizes these three neurodynamic optimization methods and their target CNLPs. Projection \eqref{eq: projection 1} is presented in the original paper of \cite{xia2007new}, and we add the Projection \eqref{eq: linear constraint projection} to facilitate the use of deep learning later. All three neurodynamic methods have the global convergence property. We refer the reader to \citep{xia2007new, qin2014two, xu2020neurodynamic} for the proof of the global convergence theorems and other details.

\subsection{Initial value problem}\label{section: 2.3}
% IVP
\DW{An initial value problem (IVP) is an ODE system together with an initial point and a time range.
The solution to the IVP is called a state solution that satisfies the initial point and the ODE system over the time range.}

% Numerical method
\DW{Almost all the ODE systems considered in this paper are nonlinear and cannot be solved analytically. 
Therefore, in practice, the IVP is usually solved by numerical integration methods, which approximate the state solution by the discretization of the domain. \citep{butcher2016numerical}.
As a typical example, Runge-Kutta methods numerically integrate the ODE system by starting with the initial point and moving forward until the desired final time is reached. 
The numerical integration method chooses a number of time points in the domain, called collocation points, and then find a solution that satisfies the ODE system at these points.
However, these conventional methods are inefficient if only the state at the end is of interest. This is due to the significant computational work required to determine all the ahead collocation points.}

% Softwares
\DW{Numerical integration methods are divided into two categories: explicit and implicit methods.
The explicit methods determine the system's state at a later time based on the current state, e.g., RK45, RK23, and DOP853 \citep{RK45, RK23, DOP853}. The implicit methods find the solution by solving equations involving the current and later states, e.g., Radau and BDF \citep{Radau, BDF}. In addition, LSODA can switch automatically between stiff and nonstiff methods \citep{LSODA}.
Scipy provides software implementations of these methods to facilitate their use \citep{2020SciPy-NMeth}.}

\section{OINN model}\label{section: 3}

\begin{figure}[t]
\centering
\includegraphics{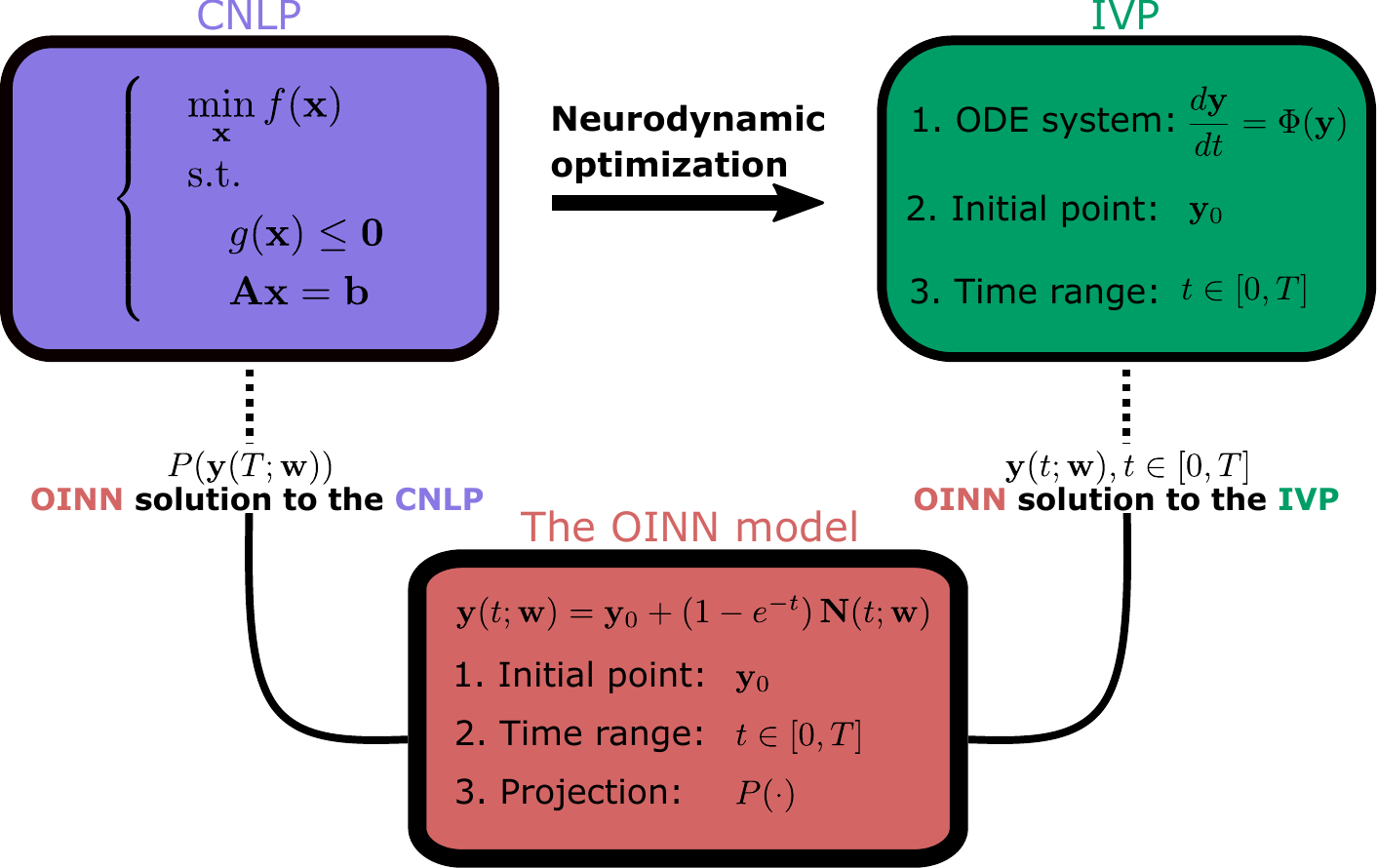}
\caption{\DW{\textbf{Problem set-up and OINN solution}
The demonstrated CNLP is a standard CNLP, where $\mathbf{x}$ is primal variables, and $\mathbf{u}$ is dual variables.
$\mathbf{y}$ is composed of $\mathbf{x}$ and $\mathbf{u}$, i.e.,  $\mathbf{y}=[\mathbf{x}, \mathbf{u}]$.}}
\label{fig: 3.1}
\end{figure}

% NLP-ODE-IVP-OINN
\DW{OINN is a generic framework for solving different CNLPs by working with neurodynamic optimization methods.
As shown in Figure\ref{fig: 3.1}, a standard CNLP is first reformulated as an IVP.
Then, an OINN model is built to solve both the CNLP and IVP. 
Let $\mathbf{y}^*$ be the optimal solution of the CNLP, $\Bar{\mathbf{y}}(t)$ be the state solution of the IVP.
In this section, we show how an OINN model provides approximations for $\mathbf{y}^*$ and $\Bar{\mathbf{y}}(t)$.}

\DW{\textbf{OINN solution to the IVP}} \ The OINN model is defined as follows
\begin{equation}
\mathbf{y}\left(t ; \mathbf{w}\right)=\mathbf{y}_0 + (1-e^{-t}) \mathbf{N}\left(t ; \mathbf{w}\right),
\label{eq: y}
\end{equation}
where $\mathbf{y}_0\in\mathbb{R}^n$ is an initial point, $[0, T]\subseteq \mathbb{R}$ is a time range, and $t\in [0, T]$ is the time variable. 
$(1-e^{-t})$ ensures that the OINN model always satisfies  the initial condition, i.e.,  $\mathbf{y}(0; \mathbf{w}) = \mathbf{y}_0$. \DW{This construction method is initially introduced by \cite{712178}, and \cite{mattheakis2022hamiltonian} demonstrates that the exponential form can result in better convergence.
$\mathbf{N}(t; \mathbf{w})$ is an neural network with learnable parameters $\mathbf{w}$. 
This paper considers the fully connected network only; other network structures are worth investigating in future research.
The OINN model itself is an approximate state solution to the IVP, i.e., 
\begin{equation}
\mathbf{y}\left(t ; \mathbf{w}\right) \approx \Bar{\mathbf{y}}(t), \quad t\in[0, T]
\label{eq: OINN model}
\end{equation}}

\DW{\textbf{OINN solution to the CNLP} \  By the neurodynamic optimization method, the endpoint of the state solution is an approximation of the optimal solution of the CNLP, i.e., 
\begin{equation}
    \overline{\mathbf{y}}(T)\approx \mathbf{y}^*
\label{eq: endpoint}
\end{equation}
Combing the endpoint $t=T$ of  \eqref{eq: OINN model} and \eqref{eq: endpoint}, we have
\begin{equation}
    P\left( \mathbf{y}\left(T; \mathbf{w}\right) \right) \approx \mathbf{y}^*,
\label{eq: approximate solution}
\end{equation}
where $P(\cdot)$ is a projection function that projects the endpoint onto a feasible set, e.g., the box-constraints \eqref{eq: projection 1} and the equality-constraints \eqref{eq: linear constraint projection}. The expression  \eqref{eq: approximate solution} represents that the endpoint of the OINN model, together with a projection function, is an approximate solution to the optimal solution of the CNLP.}

\DW{Here, we discuss two newly introduced hyperparameters in OINN, namely the initial point and the time range.}

\textbf{Initial point $\mathbf{y}_0$} \ 
Any initial point can converge to the optimal solution as long as the time goes to infinity, according to the global convergence property.
Therefore, the choice of the initial point does not affect the convergence.
However, the initial point selection has a significant impact on convergence speed; the closer the initial point is to the optimal solution, the faster the state solution approaches it.

\textbf{Time range $[0, T]$} \ 
The time range determines the training difficulty and the upper limit of accuracy.
For training difficulty, since the time range is exactly the input space of the OINN model, its span determines how large an input space the OINN model needs to be trained on.
At the same time, the time range determines the location of the state solution endpoint $\bar{\mathbf{y}}(T)$, which in turn represents an upper limit of accuracy.
Therefore, the choice of the time range span is a trade-off. 
On the one hand, the long span enables the OINN model to provide a better solution, but more training iterations are necessary to achieve it.
On the other hand, the short span is simpler to train, but the OINN model might not achieve the desired accuracy, no matter how many training iterations.

\section{OINN training}\label{section: 4}

\subsection{Loss function}

\begin{figure}[!t]
    \centering
    \includegraphics{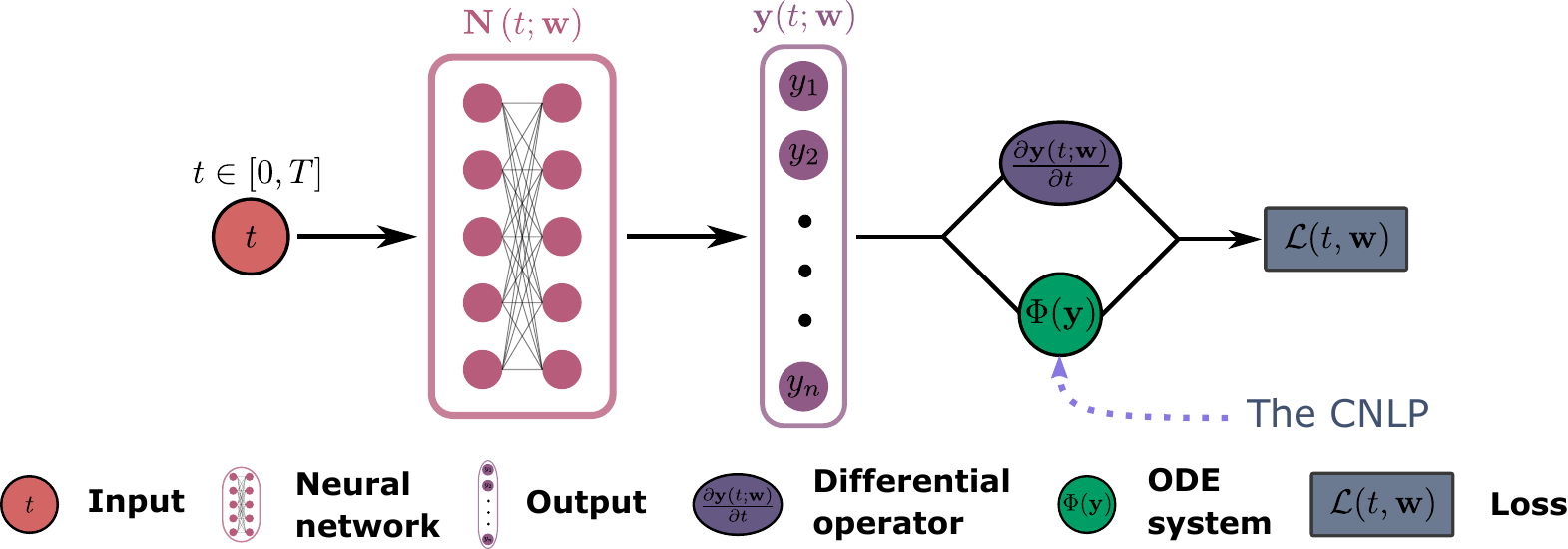}
    \caption{\DW{\textbf{Computational flow of the loss function}}}
    \label{fig: 4.1}
\end{figure}

The loss function of the OINN model is defined as follows
\begin{equation}
\mathcal{L}(t, \mathbf{w})=
e^{-\gamma*t}
\left\lVert \frac{\partial \mathbf{y}(t ; \mathbf{w})}{\partial t}
-
\Phi(\mathbf{y}(t ; \mathbf{w})) \right\rVert,
\label{eq: loss}
\end{equation}
where $\Phi(\cdot)$ refers to the ODE system corresponding to the CNLP.
\DW{$\frac{\partial \mathbf{y}(t ; \mathbf{w})}{\partial t}$ is the derivative of the output $\mathbf{y}(t; \mathbf{w})$ with respect to the input time $t$, which can be computed analytically.}
\DW{The multiplier $e^{-\gamma*t}$ reassigns the weights in the loss and gives a higher weight to the time closer to the origin, where $\gamma$ is a weighting hyperparameter. 
Using such a multiplier comes from the fact that global error can grow exponentially as a result of an early local error} \citep{flamant2020solving}. 
Figure \ref{fig: 4.1} illustrates the computational flow from a time $t$ to the loss value $\mathcal{L}(t, \mathbf{w})$.
% The loss $\mathcal{L}(t, \mathbf{w})$ measure how well the OINN model satisfies the ODE system $\Phi(\cdot)$ at time $t$ with current model parameters $\mathbf{w}$.

The objective function for the OINN model is given by
\begin{equation}
E(\mathbf{w}) =  \int_{0}^{T} \mathcal{L}(t, \mathbf{w}) d t.
\label{eq: E}
\end{equation}
The objective function $E(\mathbf{w})$ is an integral of the loss function over the time range $[0, T]$.
% By integration, the objective function eliminates the time variable $t$ and focuses only on the model parameter $\mathbf{w}$. 
The loss value $\mathcal{L}(t, \mathbf{w})$ denotes the error of the OINN model at the time $t$. The $E(\mathbf{w})$ denotes the overall error of the OINN model over the time range $[0, T]$.

\DW{However, training the model by applying gradient descent on $E(\mathbf{w})$ is infeasible since the integral over $[0, T]$ is computationally intractable. One can instead train by minimizing the batch loss
\begin{equation}
    \mathcal{L}(\mathbb{T}, \mathbf{w}) = \frac{1}{|\mathbb{T}|} \sum_{t\in\mathbb{T}} \mathcal{L}(t, \mathbf{w}),
\label{eq: batch loss}
\end{equation}
where $\mathbb{T}$ is a set of time $t$ randomly drawn from $[0, T]$. $|\mathbb{T}|$ denotes to the size of the set. }
%Such that, $\mathcal{L}(\mathbb{T}, \mathbf{w})$ is an unbiased estimate of $E(\mathbf{w})$

\subsection{Epsilon metric}
\DW{We propose a method to evaluate how well the OINN solution solves the CNLP, called the epsilon metric. The epsilon metric can be defined in two different ways, depending on the particular CNLP.}

\textbf{Epsilon: Nonlinear projection equation error} \ 
For a CNLP that can be reformulated as an NPE \eqref{eq: NPE}, such as variational inequality and nonlinear complementarity problem, the epsilon value is defined as follows
\begin{equation}
    \epsilon_1(\mathbf{y}) = \lvert P_{\Omega}(\mathbf{y}-\alpha G(\mathbf{y})) - \mathbf{y} \rvert.
\label{eq: epsilon 1}
\end{equation}
The epsilon value $\epsilon_1(\mathbf{y})$ indicates how well a solution $\mathbf{y}$ satisfies the equation \eqref{eq: NPE}. 

\textbf{Epsilon: Objective value} \ 
For the standard CNLP \eqref{eq: CNLP}, the epsilon value can be defined as follows
\begin{equation}
    \epsilon_2(\mathbf{y}) = \begin{cases} f(\mathbf{x}) & \text { if } \mathbf{x} \in \mathcal{X},  \\ 
                             +\infty       & \text { otherwise },
                \end{cases}
\label{eq: epsilon 2}
\end{equation}
where $\mathbf{y} = [\mathbf{x}, \mathbf{u}]$, $\mathcal{X}$ denotes the feasible set of the standard CNLP.
When $\mathbf{x}$ is within the feasible set, the epsilon value $\epsilon_2(\mathbf{y})$ is the objective value; otherwise, it is set to $+\infty $.
By utilizing a projection function that maps $\mathbf{x}$ onto some basic feasible set, such as $P_{eq}(\mathbf{x})$ for projecting to the equality constraint set, $\epsilon_2(\mathbf{y})$ \AL{can more likely be finite.}%has a better chance of avoiding $+\infty$.

\subsection{Training algorithm}

\begin{algorithm}[ht]
\SetAlgoLined
\SetKwInOut{Hyperparameters}{Hyperparameters}
\SetKwInOut{Input}{Input}
\SetKwInOut{Output}{Output}
\Hyperparameters{An initial point $\mathbf{y}_0$, A time range $[0, T]$}
\Input{A CNLP}
\Output{The OINN model after training}
\SetKwProg{Fn}{Function}{:}{end}
\Fn{Main}{
    Derive the ODE system $\Phi(\cdot)$ corresponding to the CNLP by a neurodynamic optimization method. \\
    Initialize an OINN model $\mathbf{y}(t; \mathbf{w})$.                                   \\
    Initialize $\epsilon_\text{best}=P\left( \mathbf{y}\left(T; \mathbf{w}\right) \right)$.\\
    \While{iter $\leq$ Max iteration}{
        $\mathbb{T} \sim  U(0, T)$: \quad Uniformly sample a batch of $t$ from the interval $[0, T]$.    \\
        Forward propagation: Compute the batch loss $\mathcal{L}(\mathbb{T}, \mathbf{w})$.            \\
        Backward propagation: Update $\mathbf{w}$ by $\nabla_\mathbf{w} \mathcal{L}(\mathbb{T},\mathbf{w})$.                            \\
        Compute the epsilon value: $\epsilon_\text{temp}=P\left( \mathbf{y}\left(T; \mathbf{w}\right) \right)$.\\
        \uIf{$\epsilon_\text{temp}<\epsilon_\text{best}$}{
            $\epsilon_\text{best} = \epsilon_\text{temp}$                               \\
            Save the OINN model with parameters $\mathbf{w}$
        }
    }
}
\caption{\DW{Training of an OINN model for solving a CNLP}}
\label{Alg: 1}
\end{algorithm}

\DW{Algorithm \ref{Alg: 1} presents an optimization procedure for the objective function $E(\mathbf{w})$ \AL{combined with} the epsilon metric.
\AL{At} each iteration, a batch of $t$ is sampled uniformly from the time range $[0, T]$ as an input dataset. 
Then, the OINN model performs gradient descent on the batch loss $\mathcal{L}(\mathbb{T}, \mathbf{w})$, which is an unbiased estimate of $E(\mathbf{w})$.}

\DW{After each round of training, the epsilon value of the OINN solution to the CNLP is computed, i.e., $\epsilon_{\text{temp}} = \epsilon\left(P(\mathbf{y}(T; \mathbf{w}))\right)$, where $\epsilon(\cdot)$ is defined by either \eqref{eq: epsilon 1} or \eqref{eq: epsilon 2}. 
Throughout the training process, the algorithm maintains the lowest epsilon value, namely $\epsilon_{\text{best}}$, representing the best solution, and the corresponding model parameter is saved.
This idea is similar to the \textit{early-stopping} in deep learning, except here, we consider the epsilon value rather than the loss.}

\begin{figure}[!t]
    \centering
    \includegraphics{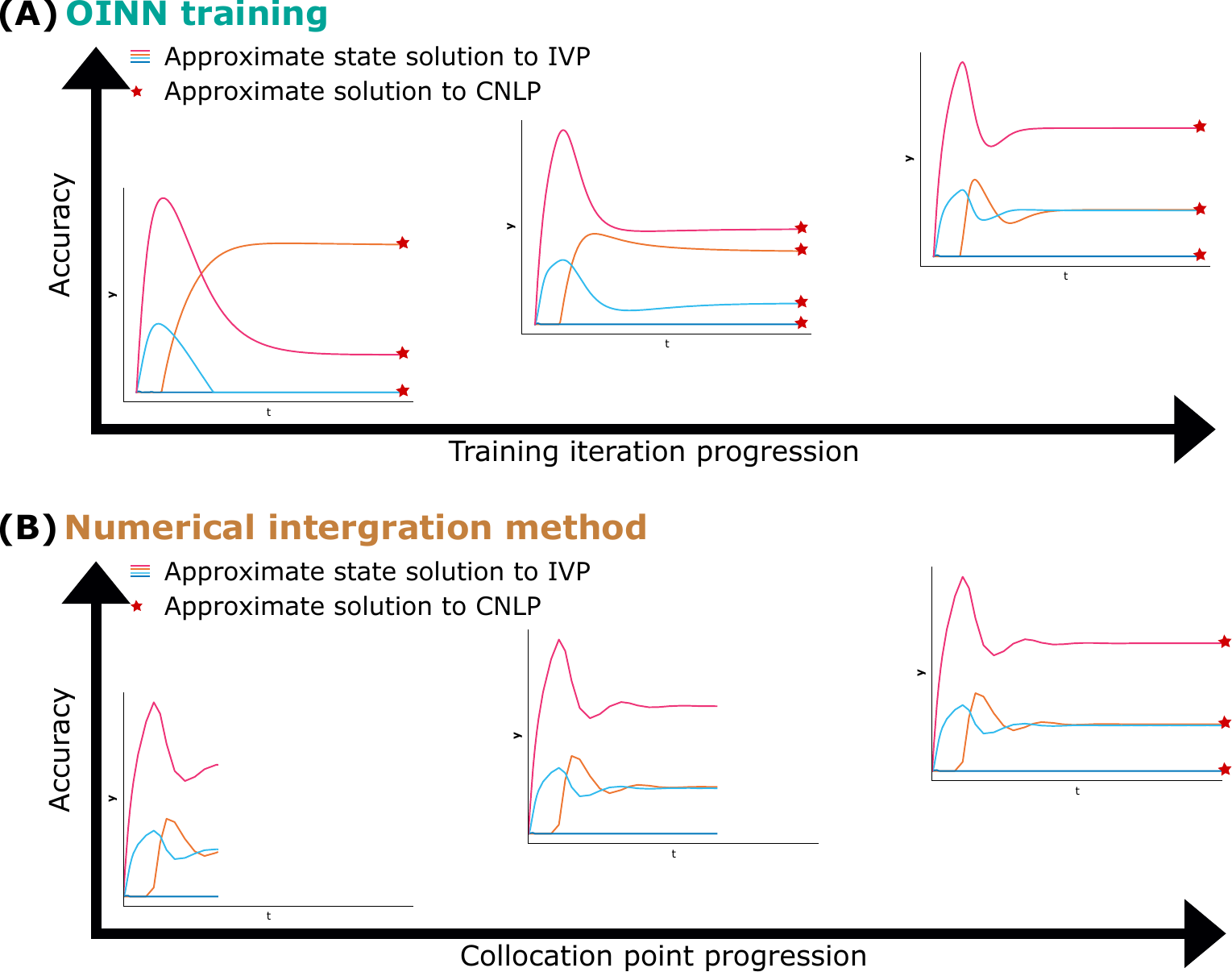}
    \caption{\DW{\textbf{Comparison between OINN and the numerical integration method for solving a CNLP}}}
    \label{fig: 2}
\end{figure}

\DW{As the OINN training progresses, the model increases its accuracy for the IVP; the prediction accuracy to the CNLP is improved by solving the IVP, as shown in Figure \ref{fig: 2}-(A).
The numerical integration method solves the IVP by stepwise integrating the ODE system and returns the solution to the CNLP at the end of the program, as shown in Figure \ref{fig: 2}-(B).
One of the promising features of OINN is that it can provide approximations for the IVP and the CNLP at any training iteration, while the numerical method can only produce solutions at the end of the program.}

\section{Numerical experiments}\label{section: 5}
% Environments
We use the Google Colab Pro+ platform to conduct our experiments, Pytorch 1.9.1 as the deep learning library. Jax 0.3.0 is used as an automatic differentiation tool  to compute the gradient or the Jacobian of a given function and subsequently models the ODE system \citep{jax2018github}.

% Hyperparameters 
The training hyperparameters are as follows
\begin{itemize}
\item The optimizer is ADAM with a learning rate of $0.001$. The decay weighting is $\gamma=0.5$.
\item The batch size is 512, and the maximum number of iterations is 50000.
\item The structure of each OINN model is a fully-connected neural network with one hidden layer of 100 neurons and Tanh as the activation function.
\end{itemize}

% Subsections
\AL{In the following, }\DW{Section \ref{section: 5.1} \AL{shows} how to use OINN to solve six different CNLP examples. Section \ref{section: 5.2} performs a hyperparameter study on the initial point $\mathbf{y}_0$ and time range $[0, T]$.
Section \ref{section: 5.3} compares the OINN method with the numerical integration methods.}

\subsection{Six CNLP examples}\label{section: 5.1}
\subsubsection{Quadratic programming}\label{section: 5.1.1}

\DW{\textbf{Example 1.} \  Consider the following quadratic programming problem}

\DW{\begin{equation}
\begin{aligned}&
\min\limits_\mathbf{x} f(\mathbf{x})= \frac{1}{2}\mathbf{x}^TQ\mathbf{x} + p^T\mathbf{x}
 \\
&\text{s.t.}           \\
&\quad  C\mathbf{x} \leq d    \\
&\quad  \mathbf{x} \geq 0,
\end{aligned}
\label{eq: Example 1}
\end{equation}
where 
\begin{equation*}
Q=\left[\begin{array}{ccc}
18 & 9 & 13 \\
9 & 14 & 6 \\
13 & 6 & 10
\end{array}\right], \quad 
p=\left[\begin{array}{c}
-30 \\
-30 \\
15
\end{array}\right], \quad 
C=\left[\begin{array}{ccc}
4 & -5 & -4\\
-5 & -2 & -4
\end{array}\right], \quad 
d=\left[\begin{array}{c}
-5 \\
1 
\end{array}\right].
\end{equation*}}

\DW{The CNLP \eqref{eq: Example 1} can be reformulated as the NPE as follows
\begin{equation}
\left(\mathbf{y} - \left(M\mathbf{y}+q\right)\right)^+=\mathbf{y},
\label{eq: QP-NPE}
\end{equation}
where $\mathbf{y}=[x_1, x_2, x_3, u_1, u_2]^T$, $x_1$, $x_2$ and $x_3$ are decision variables, and $u_1$, $u_2$ are dual variables. $(\mathbf{y})^+ = \max\{\mathbf{0}, \mathbf{y}\}$. $M$ and $q$ are denoted as
\begin{equation*}
M=\left[\begin{array}{cc}
Q & C^{\mathrm{T}} \\
-C & 0
\end{array}\right], \quad 
q=\left[\begin{array}{l}
p \\
d
\end{array}\right].
\end{equation*}}

\DW{The following ODE system model the NPE
\begin{equation}
\frac{\mathrm{d} \mathbf{y}}{\mathrm{d} t}= - M(\mathbf{y})^+ - q +(\mathbf{y})^+-\mathbf{y}.
\label{eq: ODE-Example1}
\end{equation}
The ODE system together with the initial point $\mathbf{y}_ 0=[0, 0, 0, 0, 0]$ and time range $[0, 10]$ form the IVP as follow
\begin{equation}
\eqref{eq: ODE-Example1}, \quad \mathbf{y}_ 0=[0, 0, 0, 0, 0], \quad t\in[0, 10]
\label{eq: IVP-Example1}
\end{equation}}

\DW{An OINN model, $\mathbf{y}(t; \mathbf{w}) \ t\in [0, 10]$, is built as an approximate state solution to the IVP \eqref{eq: IVP-Example1}, and its endpoint $\left(\mathbf{y}(10; \mathbf{w})\right)^+$ is an approximate solution to the NPE \eqref{eq: QP-NPE}.}

\begin{figure}[h]
    \centering
    \includegraphics{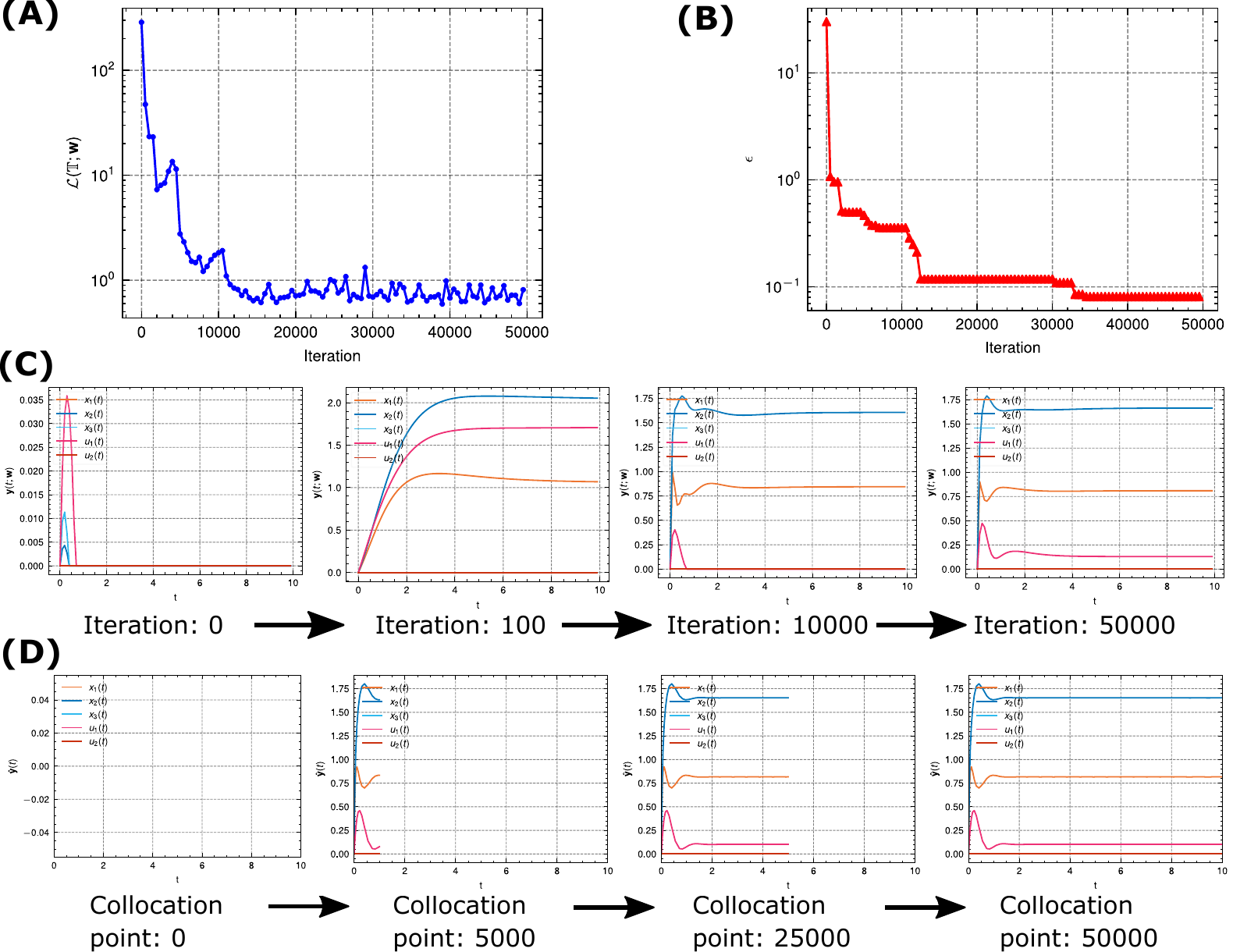}
    \caption{\textbf{Example 1 Quadratic programming} 
    \textbf{(A) The loss versus the number of iterations.} $\mathcal{L}\left(\mathbb{T}, \mathbf{w}\right)$ refers to the batch loss defined in \eqref{eq: batch loss} and \eqref{eq: loss} 
    \textbf{(B) The epsilon value versus the number of iterations.} The epsilon metric is defined in \eqref{eq: epsilon 1}
    \textbf{(C) The solving process of the OINN model}
    \textbf{(D) The solving process of the  numerical integration method}}
    \label{fig: Example 1}
\end{figure}

\begin{table}[h]
\begin{tabular}{@{}l|ll|ll@{}}
\toprule
\multirow{2}{*}{Index}     & \multicolumn{2}{l|}{OINN}                       & \multicolumn{2}{l}{Numerical integration method}      \\ \cmidrule(l){2-5} 
                           & Iteration & Solution                            & Collocation point & Solution \\ \midrule
\multirow{6}{*}{Example 1} & 0         & [0.00, 0.00, 0.00, 0.00, 0.00]     & 0                  & [0.00, 0.00, 0.00, 0.00, 0.00]         \\
                           & 10        & [0.00, 0.42, 0.00, 0.17, 0.00]     & 10                 & [0.06, 0.06, 0.00, 0.01, 0.00]         \\
                           & 100       & [1.07, 2.05, 0.00, 1.71, 0.00]     & 100                &  [0.46, 0.49, 0.00, 0.09, 0.00]        \\
                           & 1000      & [0.75, 1.73, 0.00, 0.00, 0.00]     & 1000               &  [0.82, 1.63, 0.00, 0.46, 0.00]        \\
                           & 10000     & [0.84, 1.61, 0.00, 0.00, 0.00]     & 10000              &  [0.82, 1.65, 0.00,  0.10, 0.00]        \\
                           & 50000     & [0.81, 1.66, 0.00, 0.13,  0.00]     & 50000             &  [0.82, 1.65, 0.00, 0.10, 0.00]        \\ \bottomrule
\end{tabular}
\caption{
\textbf{Example 1, Approximate solutions to the NPE during solving}\ 
We choose a step size of 0.0002 for the numerical integration method. collocation points 0, 10, 100, 1000, 10000, 50000 represent the time ranges $[0, 0]$, $[0, 0.002]$, $[0, 02]$, $[0, 2]$, $[0, 10]$ respectively.}
\label{tab: Example 1}
\end{table}

\DW{Figure \ref{fig: Example 1} shows the training of this OINN model,  where the loss decreased from the initial value of 287.26 to 0.62, and the epsilon value decreased from 30.00 to 0.08.
Figure \ref{fig: Example 1} (C) and (D) show the progressions of the approximate state solutions to the IVP \eqref{eq: IVP-Example1}.}

\DW{Table \ref{tab: Example 1} displays the progressions of the approximate solutions to the NPE \eqref{eq: QP-NPE}. 
The OINN model gives the final solution of $[0.81, 1.66, 0.00, 0.13, 0.00]$ to the NPE, where $[0.81, 1.66, 0.00]$ is the solution to the CNLP \eqref{eq: Example 1}.
The numerical integration method gives the final solution of $[0.82, 1.65, 0.00, 0.10, 0.00]$ to the NPE, where $[0.82, 1.65, 0.00]$ is the solution to the CNLP \eqref{eq: Example 1}.}

\subsubsection{Convex-smooth standard CNLP}\label{section: 5.1.2}

\DW{\textbf{Example 2.} Consider the following convex-smooth standard CNLP:
\begin{equation}
\begin{aligned}
&\min\limits_\mathbf{x} f(\mathbf{x}) = x_1^2+2x_2^2+2x_1x_2-10x_1-12x_2 \\
&\text{s.t.}           \\
&\quad g_1(\mathbf{x}) = x_1+3x_2-8 \leq 0    \\
&\quad g_2(\mathbf{x}) = x_1^2+x_2^2+2x_1-2x_2-3 \leq 0 \\
&\quad  0 \leq \mathbf{x} \leq 2.
\end{aligned}    
\label{eq: Example 2}
\end{equation}}

\DW{The CNLP can be reformulated as the following NPE
\begin{equation}
P_{\Omega}(\mathbf{y} - G(\mathbf{y}))=\mathbf{y},
\label{eq: Example2-NPE}
\end{equation}
where $\mathbf{y}=[x_1, x_2, u_1, u_2]^T$; $x_1$, $x_2$ are decision variables, and $u_1$, $u_2$  are dual variables. $P_{\Omega}(\mathbf{y})$ is a projection function defined in \eqref{eq: projection 1} which projects $\mathbf{y}\in \mathbb{R}^4$ onto the set $\Omega=\{ \mathbf{y}\in \mathbb{R}^4 \mid 0\leq x_1 \leq 2, \ 0\leq x_2 \leq 2, \ u_1\geq0, \  u_2\geq0\}$. $G(\mathbf{y})$ is defined as
\begin{equation}
G(\mathbf{y})=\left[\begin{array}{l}
\nabla f(\mathbf{x}) + \nabla g(\mathbf{x})^T\mathbf{u} \\
-g(\mathbf{x})
\end{array}\right],
\end{equation}
where $g(\mathbf{x})=[g_1(\mathbf{x}), g_2(\mathbf{x})]^T$, $\mathbf{x}=[x_1, x_2]^T$, $\mathbf{u}=[u_1, u_2]^T$.}

\DW{The following ODE system \AL{models} the NPE \eqref{eq: Example2-NPE}
\begin{equation}
\frac{\mathrm{d} \mathbf{y}}{\mathrm{d} t}= - G\left(P_{\Omega}(\mathbf{y})\right)+P_{\Omega}(\mathbf{y})-\mathbf{y},
\label{eq: ODE-Example2}
\end{equation}
The ODE system together with the initial point $\mathbf{y}_ 0=[0, 0, 0, 0]$ and time range $[0, 10]$ form an IVP as follow
\begin{equation}
\eqref{eq: ODE-Example2}, \quad \mathbf{y}_ 0=[0, 0, 0, 0], \quad t\in[0, 10]
\label{eq: IVP-Example2}
\end{equation}}

\DW{An OINN model, $\mathbf{y}(t; \mathbf{w}) \ t\in [0, 10]$, is built as an approximate state solution to this IVP \eqref{eq: IVP-Example2}, and its endpoint $P_\Omega\left(\mathbf{y}(10; \mathbf{w})\right)$ is an approximate solution to the NPE \eqref{eq: Example2-NPE}.}

\begin{figure}[!t]
    \centering
    \includegraphics{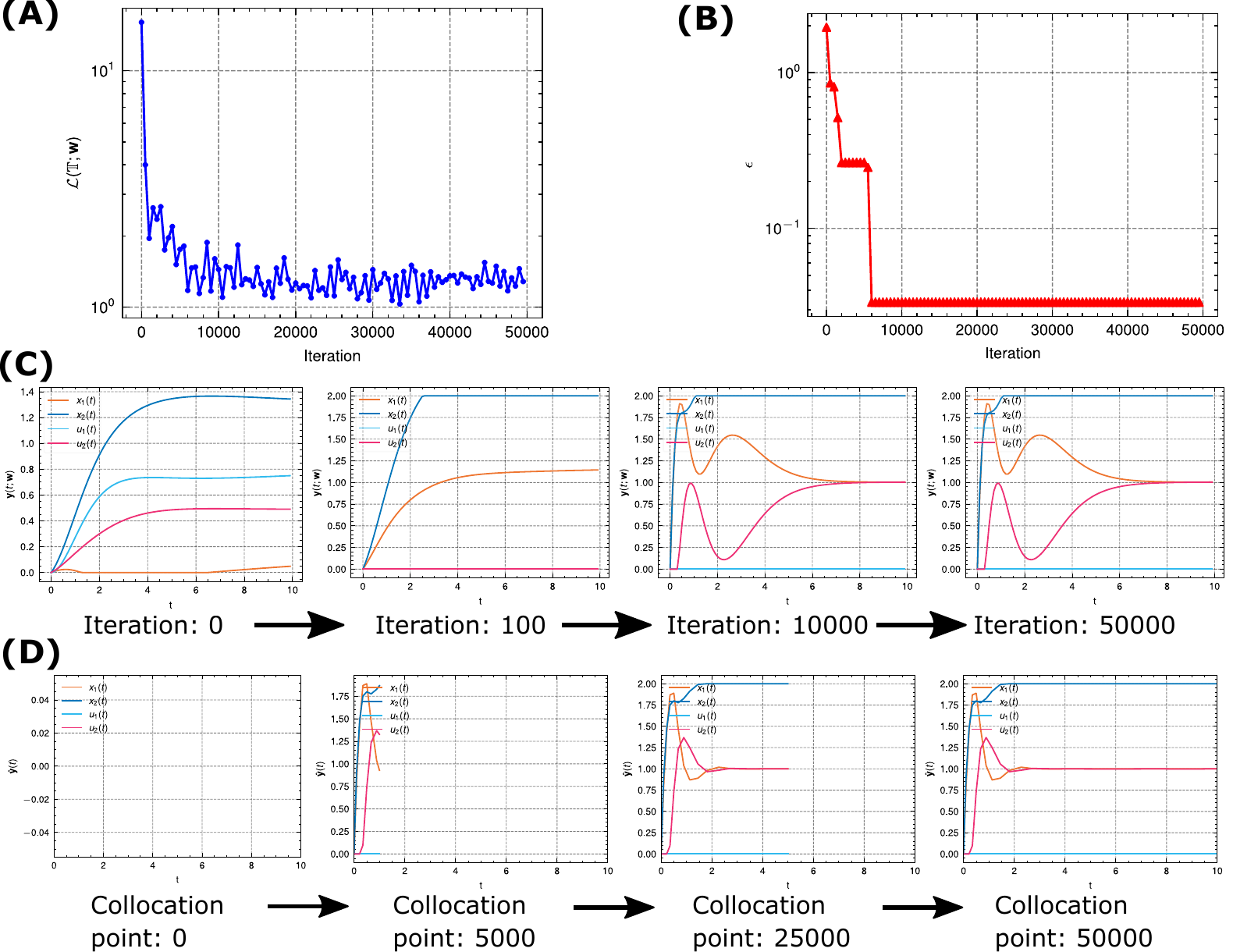}
    \caption{\textbf{Example 2 Convex smooth standard CNLP} 
    \textbf{(A) The loss versus the number of iterations.} 
    \textbf{(B) The epsilon value versus the number of iterations.} The epsilon metric is defined in \eqref{eq: epsilon 1}
    \textbf{(C) The solving process of the OINN model}
    \textbf{(D) The solving process of the  numerical integration method}}
    \label{fig: Example 2}
\end{figure}

\begin{table}[!ht]
\begin{tabular}{@{}l|ll|ll@{}}
\toprule
\multirow{2}{*}{Index}     & \multicolumn{2}{l|}{OINN}                       & \multicolumn{2}{l}{Numerical integration method}      \\ \cmidrule(l){2-5} 
                           & Iteration & Solution                            & Collocation point  & Solution \\ \midrule
\multirow{6}{*}{Example 2} & 0         & [0.05, 1.34, 0.75, 0.49]            & 0                  & [0.00, 0.00, 0.00, 0.00]         \\
                           & 10        & [0.84, 2.00, 0.00, 0.00]            & 10                 & [0.02, 0.02, 0.00, 0.00]         \\
                           & 100       & [1.15, 2.00, 0.00, 0.00]            & 100                &  [0.19, 0.23, 0.00, 0.00]        \\
                           & 1000      & [1.19, 2.00, 0.00, 0.00]            & 1000               &  [1.36, 1.42, 0.00, 0.00]        \\
                           & 10000     & [1.00, 2.00, 0.00, 1.00]            & 10000              &  [1.01, 2.00, 0.00, 0.97]        \\
                           & 50000     & [1.00, 2.00, 0.00, 1.00]            & 50000              &  [1.00, 2.00, 0.00, 1.00]        \\ \bottomrule
\end{tabular}
\centering
\caption{\textbf{Example 2, Approximate solutions to the NPE during solving}}
\label{tab: Example 2}
\end{table}

\DW{Figure  \ref{fig: Example 2} shows the training of the OINN model, where the loss value decreased from $16.02$ to $1.11$, and the epsilon value decreased from $1.95$ to $0.03$.
Figure \ref{fig: Example 2} (C) and (D) show the progressions of the approximate state solutions to the IVP \eqref{eq: IVP-Example2}.}

\DW{Table \ref{tab: Example 2} displays the progression of the approximate solutions to the NPE \eqref{eq: Example2-NPE}. 
Both the OINN model and the numerical integration method gives the same final solution of $[1.00, 2.00, 0.00, 1.00]$ to the NPE \eqref{eq: Example2-NPE}, where $[1.00, 2.00]$ is the solution to the CNLP \eqref{eq: Example 2}.}

\subsubsection{Variational inequality}\label{section: 5.1.3}

\DW{\textbf{Example 3.} Consider the following variational inequality
\begin{equation}
    \left(\mathbf{y}-\mathbf{y}^{*}\right)^{T} G\left(\mathbf{y}^{*}\right) \geq 0, \quad \mathbf{y} \in \Omega,
\label{eq: Example 3}
\end{equation}
where 
\begin{equation*}
    G(\mathbf{y}) = \left[\begin{array}{c}
y_1 - \frac{2}{(y_1+0.8)} + 5y_2 - 13\\
1.2y_1 + 7y_2\\
3y_3 + 8y_4 \\
y_3 + 2y_4 - \frac{4}{(y_4+2)} - 12
\end{array}\right], \quad
\begin{aligned}\Omega = \{\mathbf{y} \in \mathbb{R}^{4} \mid & 1\leq y_1\leq100, -3\leq y_2\leq100, \\
                                                 &-3\leq y_3\leq100, 1\leq y_4\leq100 \}.
  \end{aligned}
\end{equation*}.}

\DW{The problem can be reformulated as the following NPE
\begin{equation}
P_{\Omega}(\mathbf{y} - G(\mathbf{y}))=\mathbf{y}.
\label{eq: Example3-NPE}
\end{equation}}

\DW{The following ODE system model the NPE \eqref{eq: Example3-NPE}
\begin{equation}
\frac{\mathrm{d} \mathbf{y}}{\mathrm{d} t}= - G\left(P_{\Omega}(\mathbf{y})\right)+P_{\Omega}(\mathbf{y})-\mathbf{y}.
\label{eq: ODE-Example3}
\end{equation}
The ODE system together with the initial point $\mathbf{y}_ 0=[0, 0, 0, 0]$ and time range $[0, 10]$ form the IVP as follow
\begin{equation}
\eqref{eq: ODE-Example3}, \quad \mathbf{y}_ 0=[0, 0, 0, 0], \quad t\in[0, 10].
\label{eq: IVP-Example3}
\end{equation}}

\DW{An OINN model, $\mathbf{y}(t; \mathbf{w}) \ t\in [0, 10]$, is built as an approximate state solution to this IVP \eqref{eq: IVP-Example3}, and its endpoint $P_\Omega\left(\mathbf{y}(10; \mathbf{w})\right)$ is an approximate solution to the NPE \eqref{eq: Example3-NPE}.}

\begin{figure}[!ht]
    \centering
    \includegraphics{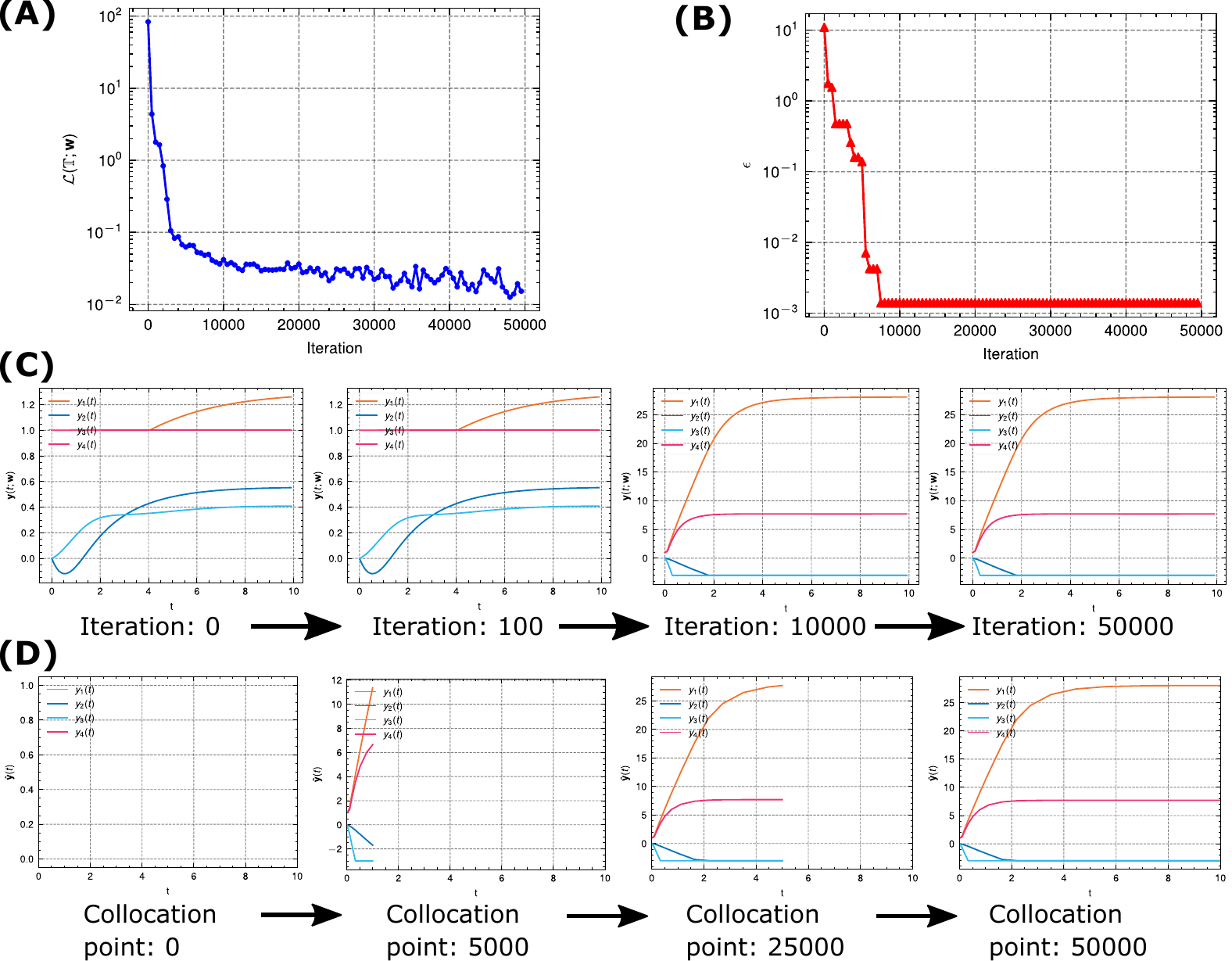}
    \caption{\textbf{Example 3 variational inequality} 
    \textbf{(A) The loss versus the number of iterations.} 
    \textbf{(B) The epsilon value versus the number of iterations.} The epsilon metric is defined in \eqref{eq: epsilon 1}
    \textbf{(C) The solving process of the OINN model}
    \textbf{(D) The solving process of the  numerical integration method}}
    \label{fig: Example 3}
\end{figure}

\begin{table}[!ht]
\begin{tabular}{@{}l|ll|ll@{}}
\toprule
\multirow{2}{*}{Index}     & \multicolumn{2}{l|}{OINN}                       & \multicolumn{2}{l}{Numerical integration method}      \\ \cmidrule(l){2-5} 
                           & Iteration & Solution                            & Collocation point & Solution \\ \midrule
\multirow{6}{*}{Example 3} & 0         & [1.26, 0.55, 0.41, 1.00]            & 0                 & [1.00, 0.00, 0.00, 1.00]         \\
                           & 10        & [1.26, 0.55, 0.41, 1.00]            & 10                & [1.00, 0.00, -0.02,  1.00]         \\
                           & 100       & [1.26, 0.55, 0.41, 1.00]            & 100               &  [1.00, -0.02, -0.16,  1.00]        \\
                           & 1000      & [26.52, -3.00, -3.00, 7.67]         & 1000              &  [2.61, -0.22, -1.74,  2.28]        \\
                           & 10000     & [28.07, -3.00, -3.00, 7.71]         & 10000             &  [20.79, -3.00, -3.00, 7.57]        \\
                           & 50000     & [28.07, -3.00, -3.00, 7.71]         & 50000             &  [28.06, -3.00, -3.00, 7.70]        \\ \bottomrule
\end{tabular}
\centering
\caption{\textbf{Example 3, Approximate solutions to the NPE during solving }}
\label{tab: Example 3}
\end{table}

\DW{Figure  \ref{fig: Example 3} shows the training of the OINN model, where the loss value decreased from $83.38$ to $0.01$, and the epsilon value decreased from $10.92$ to $0.00$. Figure \ref{fig: Example 3} (C) and (D) show the progression of the approximate state solutions to the IVP \eqref{eq: IVP-Example3}.}

\DW{Table \ref{tab: Example 3} displays the progression of the approximate solutions to the NPE \eqref{eq: Example3-NPE}. 
The OINN model gives the final solution of $[28.07, -3.00, -3.00, 7.71]$ to both the variational inequality \eqref{eq: Example 3} and NPE \eqref{eq: Example3-NPE}.
The numerical method gives the final solution of $[28.06, -3.00, -3.00, 7.70]$.}

\subsubsection{Nonlinear complementary problem}\label{section: 5.1.4}

\DW{\textbf{Example 4} Consider the following nonlinear complementary problem 
\begin{equation}
\mathbf{y}^{T} F(\mathbf{y})=0, \quad F(\mathbf{y}) \geq 0, \quad \mathbf{y} \geq 0,
\label{eq: Example 4}
\end{equation}
where 
\begin{equation*}
F(\mathbf{y})= \left(\qquad \begin{array}{c}
2 y_{1} e^{\left(y_{1}^{2}+\left(y_{2}-1\right)^{2}\right)}+y_{1}-y_{2}-y_{3}+1 \\
2\left(y_{2}-1\right) e^{\left(y_{1}^{2}+\left(y_{2}-1\right)^{2}\right)}-y_{1}+2 y_{2}+2 y_{3}+3 \\
-y_{1}+2 y_{2}+3 y_{3}
\end{array}\right).
\end{equation*}}

\DW{The problem can be reformulated as the following NPE
\begin{equation}
(\mathbf{y}-F(\mathbf{y}))^+=\mathbf{y}.
\label{eq: Example4-NPE}
\end{equation}}

\DW{The following ODE system model the NPE \eqref{eq: Example4-NPE}
\begin{equation}
\frac{\mathrm{d} \mathbf{y}}{\mathrm{d} t}= - F\left((\mathbf{y})^+\right)+(\mathbf{y})^+-\mathbf{y},
\label{eq: ODE-Example4}
\end{equation}
The ODE system together with the initial point $\mathbf{y}_ 0=[0, 0, 0]$ and time range $[0, 10]$ form the IVP as follow
\begin{equation}
\eqref{eq: ODE-Example4}, \quad \mathbf{y}_ 0=[0, 0, 0], \quad t\in[0, 10]
\label{eq: IVP-Example4}
\end{equation}}

\DW{An OINN model, $\mathbf{y}(t; \mathbf{w}) \ t\in [0, 10]$, is built as an approximate state solution to this IVP \eqref{eq: IVP-Example4}, and its endpoint $P_\Omega\left(\mathbf{y}(10; \mathbf{w})\right)$ is an approximate solution to the NPE \eqref{eq: Example4-NPE}.}

\begin{figure}[!ht]
    \centering
    \includegraphics{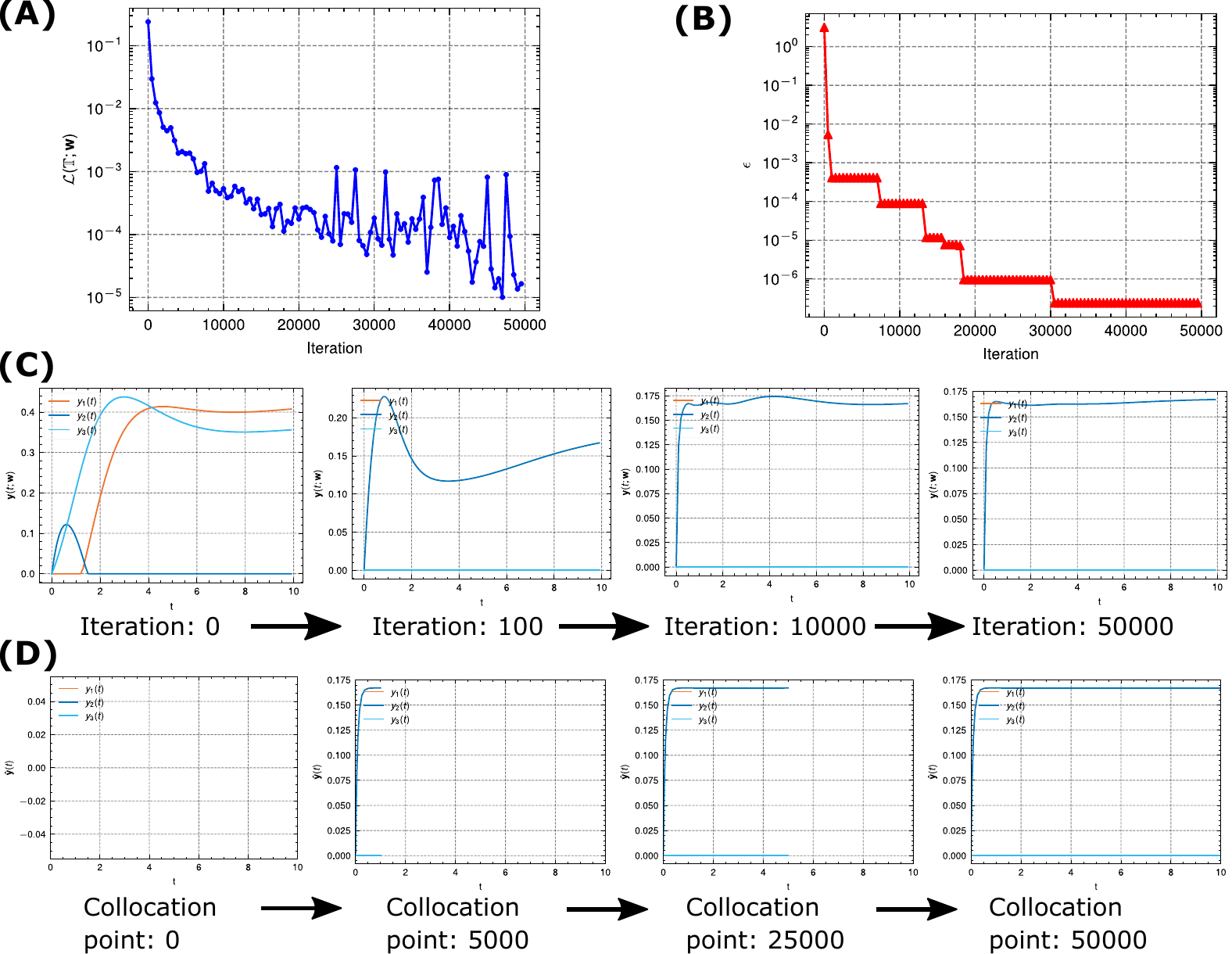}
    \caption{\textbf{Example 4 nonlinear complementary problem} 
    \textbf{(A) The loss versus the number of iterations.} 
    \textbf{(B) The epsilon value versus the number of iterations.} The epsilon metric is defined in \eqref{eq: epsilon 1}
    \textbf{(C) The solving process of the OINN model}
    \textbf{(D) The solving process of the  numerical integration method}}
    \label{fig: Example 4}
\end{figure}

\DW{Figure  \ref{fig: Example 4} shows the training of the OINN model, where the loss value decreased from $0.24$ to $0.00$, and the epsilon value decreased from $3.12$ to $0.00$. 
Figure \ref{fig: Example 4} (C) and (D) display the progression of the approximate state solutions to the IVP \eqref{eq: IVP-Example4}.}

\begin{table}[!ht]
\begin{tabular}{@{}l|ll|ll@{}}
\toprule
\multirow{2}{*}{Index}     & \multicolumn{2}{l|}{OINN}                       & \multicolumn{2}{l}{Numerical integration method}      \\ \cmidrule(l){2-5} 
                           & Iteration & Solution                            & Collocation point & Solution \\ \midrule
\multirow{6}{*}{Example 4} & 0         & [0.41, 0.00, 0.36]                  & 0                 & [0.00, 0.00, 0.00]         \\
                           & 10        & [0.00, 0.11, 0.00]                  & 10                & [0.00, 0.00, 0.00]         \\
                           & 100       & [0.00, 0.17, 0.00]                  & 100               &  [0.00, 0.04, 0.00]        \\
                           & 1000      & [0.00, 0.17, 0.00]                  & 1000              &  [0.00, 0.15, 0.00]        \\
                           & 10000     & [0.00, 0.17, 0.00]                  & 10000             &  [0.00, 0.17, 0.00]        \\
                           & 50000     & [0.00, 0.17, 0.00]                  & 50000             &  [0.00, 0.17, 0.00]        \\ \bottomrule
\end{tabular}
\centering
\caption{\textbf{Example 4, Approximate solutions to the NPE during solving}}
\label{tab: Example 4}
\end{table}

\DW{Table \ref{tab: Example 4} displays the progression of the approximate solutions to the NPE \eqref{eq: Example4-NPE}. 
The OINN model gives the final solution of $[0.00, 0.17, 0.00]$ for both the nonlinear complementary problem \eqref{eq: Example 4} and NPE \eqref{eq: Example4-NPE}.}

\subsubsection{Convex nonsmooth standard CNLP}\label{section: 5.1.5}

\DW{\textbf{Example 5} Consider the following convex nonsmooth standard CNLP
\begin{equation}
    \begin{aligned}
&\min\limits_\mathbf{x} f(\mathbf{x}) = 10(x_1+x_2)^2 + (x_1-2)^2 + 20|x_3-3| + e^{x_3} \\
&\text{s.t.}           \\
&\quad  g(\mathbf{x}) = (x_1+3)^2 + x_2  \leq 36    \\
&\quad  h(\mathbf{x}) = 2x_1+5x_3 - 7 = 0.
\end{aligned}
\label{eq: Example 5}
\end{equation}
Denote $\mathbf{y} = [x_1, x_2, x_3, u]^T$, where $x_1$, $x_2$, $x_3$ are primal variables, and $u$ is dual variable. Denote $\mathbf{A}=[2, 0, 5]$, $\mathbf{b}=7$, $\mathbf{U}=\mathbf{A}^{T}\left(\mathbf{A} \mathbf{A}^{T}\right)^{-1} \mathbf{A}$, and $\mathbf{I}_3$ is the identity matrix of size $3 \times 3$. }

\DW{The following ODE system \AL{models} this CNLP
\begin{equation}
    \begin{aligned} 
        \frac{d\mathbf{x}}{dt} = &-(\mathbf{I}_3-\mathbf{U})\left(\nabla f(\mathbf{x})+ (u+g(\mathbf{x}))^{+}\nabla g(\mathbf{x}) \right) - \mathbf{A}^{T} h(\mathbf{x}), \\ 
        \frac{du}{dt} =& \frac{1}{2}\left(-u+(u+g(\mathbf{x}))^{+}\right).
    \end{aligned}
\label{eq: ODE-Example5}
\end{equation}
The ODE system together with the initial point $\mathbf{y}_ 0=[0, 0, 0, 0]$ and time range $[0, 10]$ form the IVP as follow
\begin{equation}
\eqref{eq: ODE-Example5}, \quad \mathbf{y}_ 0=[0, 0, 0, 0], \quad t\in[0, 10]
\label{eq: IVP-Example5}
\end{equation}}

\DW{An OINN model, $\mathbf{y}(t; \mathbf{w}) \ t\in [0, 10]$, is built as an approximate state solution to this IVP, and its endpoint $P_{eq}\left(\mathbf{y}(10; \mathbf{w})\right)$ is an approximate solution to the CNLP. $P_{eq}(\cdot)$ is a projection function used to project $\mathbf{x}$ onto the equality constraint set $\{\mathbf{x} \in \mathbb{R}^3 | \ h(\mathbf{x})=0\}$, as defined in \eqref{eq: linear constraint projection}.}

\begin{figure}[!ht]
    \centering
    \includegraphics{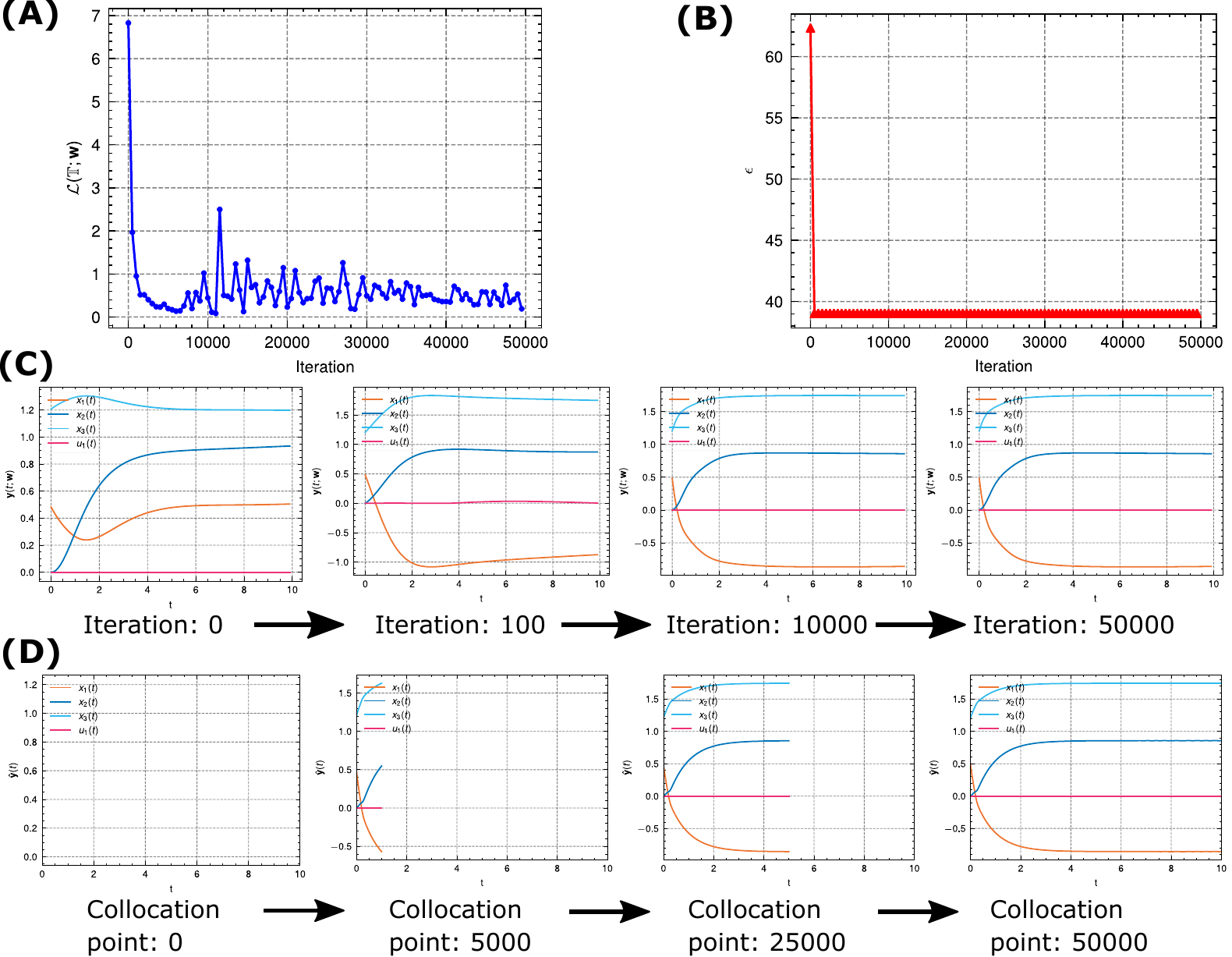}
    \caption{\textbf{Example 5 convex nonsmooth standard CNLP} 
    \textbf{(A) The loss versus the number of iterations.} 
    \textbf{(B) The epsilon value versus the number of iterations.} The epsilon metric is defined in \eqref{eq: epsilon 2}
    \textbf{(C) The solving process of the OINN model}
    \textbf{(D) The solving process of the  numerical integration method}}
    \label{fig: Example 5}
\end{figure}

\begin{table}[!ht]
\begin{tabular}{@{}l|ll|ll@{}}
\toprule
\multirow{2}{*}{Index}     & \multicolumn{2}{l|}{OINN}                       & \multicolumn{2}{l}{Numerical integration method}      \\ \cmidrule(l){2-5} 
                           & Iteration & Solution                            & Collocation point & Solution \\ \midrule
\multirow{6}{*}{Example 5} & 0         & [0.51, 0.93, 1.20,  0.00]           & 0                 & [0.48, 0.00, 1.21, 0.00]         \\
                           & 10        & [-0.30,   0.42,  1.52,  0.34]       & 10                & [0.48, 0.00,  1.21, 0.00]         \\
                           & 100       & [-0.87,  0.87,  1.75,  0.00]        & 100               &  [0.42, 0.00, 1.23, 0.00]        \\
                           & 1000      & [-0.87,  0.87,  1.75,  0.00]        & 1000              &  [0.00,  0.07, 1.40,  0.00]        \\
                           & 10000     & [-0.86,  0.86,  1.74,  0.00]        & 10000             &  [-0.78,  0.77,  1.71,  0.00]        \\
                           & 50000     & [-0.86,  0.86,  1.74,  0.00]        & 50000             &  [-0.86, 0.86,  1.74,  0.00]        \\ \bottomrule
\end{tabular}
\centering
\caption{\textbf{Example 5, Approximate solutions to the CNLP during solving}}
\label{tab: Example 5}
\end{table}

\DW{Figure  \ref{fig: Example 5} shows the training of this OINN model, where the loss decreased from $6.83$ to $0.37$, and the epsilon value decreased from $62.31$ to $39.01$. 
In this example, the epsilon value is defined as the objective value, as in \eqref{eq: epsilon 2}.
Figure \ref{fig: Example 5} (C) and (D) show the progressions of the approximate state solutions to the IVP \eqref{eq: IVP-Example5}.}

\DW{Table \ref{tab: Example 5} displays the progressions of approximate solutions to the CNLP \eqref{eq: Example 5}. 
Both the OINN model and the numerical integration method give the same solution $[-0.86, 0.86, 1.74, 0.00]$, where $[-0.86, 0.86, 1.74]$ is the solution for the primal variable $\mathbf{x}$, and $0.00$ is the solution for the  dual variable $u$.}

\subsubsection{Pseudoconvex nonsmooth standard CNLP}\label{section: 5.1.6}

\DW{\textbf{Example 6} Consider the following pseudoconvex nonsmooth standard CNLP 
\begin{equation}
\begin{aligned}
&\min\limits_\mathbf{x} f(\mathbf{x}) = \frac{x_1+x_2+e^{|x_2-1|}-40}{(x_1+x_2+x_3)^2+3} \\
&\text{s.t.}           \\
&\quad  g_1(\mathbf{x}) = -3x_1+2x_2 - 5  \leq 0   \\
&\quad  g_2(\mathbf{x}) = x_1^2 + x_2 - 3 \leq 0   \\
&\quad  h(\mathbf{x}) = x_1 + 2x_2 + x_3 - 2 = 0 
\end{aligned}
\label{eq: Example 6}
\end{equation}
Denote $\mathbf{A}=[1, 2, 1]$, $b=2$, $\mathbf{U}=\mathbf{A}^{T}\left(\mathbf{A} \mathbf{A}^{T}\right)^{-1} \mathbf{A}$. }

\DW{The following ODE system model this CNLP
\begin{equation}
    \frac{d\mathbf{x}}{d t} 
    =
    -\theta(t)(\mathbf{I}_3-\mathbf{U})\left(\mu(\mathbf{x}) \nabla f(\mathbf{x}) + \partial B(\mathbf{x}) \right) -  \text{sign}(h(\mathbf{x}))\mathbf{A}^T,
\label{eq: ODE-Example6}
\end{equation}
where $\text{sign}(\cdot)$ is the sign function. $\theta(t)$ is defined by 
\begin{equation*}
\theta(t)= \begin{cases}0, & \text { if } t \leq T_{0} \\ 1, & \text { otherwise, }\end{cases}
\end{equation*}
where $T_{0}=1+\left\|\mathbf{A} x_{0}-b\right\|_{1} / \lambda_{\min }\left(\mathbf{A} \mathbf{A}^{T}\right)$, $\lambda_{\min }\left(\mathbf{A} \mathbf{A}^{T}\right)$ represents the minimum eigen value of the matrix $\mathbf{A} \mathbf{A}^T$,
$\mathbf{x}_0=[0, 0, 0]$.
$\mu(\mathbf{x})$ is defined by
\begin{equation*}
\mu(\mathbf{x})= \begin{cases}1, & \text { if } g_1(\mathbf{x}) \leq 0 \ \&\  g_2(\mathbf{x}) \leq 0, \\ 0, & \text { otherwise. }\end{cases}
\end{equation*}
$\partial B(\mathbf{x})$ is defined by
\begin{equation*}
\partial B(\mathbf{x})
= 
\begin{cases}
0, & \text { if } g_1(\mathbf{x}) \leq 0 \ \&\  g_2(\mathbf{x}) \leq 0,\\
\nabla g_1(\mathbf{x}), & \text { if } g_1(\mathbf{x}) > 0 \ \&\  g_2(\mathbf{x}) \leq 0,\\
\nabla g_2(\mathbf{x}), & \text { if } g_1(\mathbf{x}) \leq 0 \ \&\  g_2(\mathbf{x}) > 0,\\
\nabla g_1(\mathbf{x})+\nabla g_2(\mathbf{x}), & \text { if } g_1(\mathbf{x}) > 0 \ \&\  g_2(\mathbf{x}) > 0.\\
\end{cases}
\end{equation*}}

\DW{The ODE system together with the initial point $\mathbf{y}_ 0=[0, 0, 0]$ and time range $[0, 10]$ form the IVP as follow
\begin{equation}
 \eqref{eq: ODE-Example6}, \quad \mathbf{y}_ 0=[0, 0, 0], \quad t\in[0, 10]
\label{eq: IVP-Example6}
\end{equation}}

\DW{An OINN model, $\mathbf{y}(t; \mathbf{w}) \ t\in [0, 10]$, is built as an approximate state solution for the IVP, and its endpoint $P_{eq}\left(\mathbf{y}(10; \mathbf{w})\right)$ is an approximate solution to the CNLP.}

\begin{figure}[!ht]
    \centering
    \includegraphics{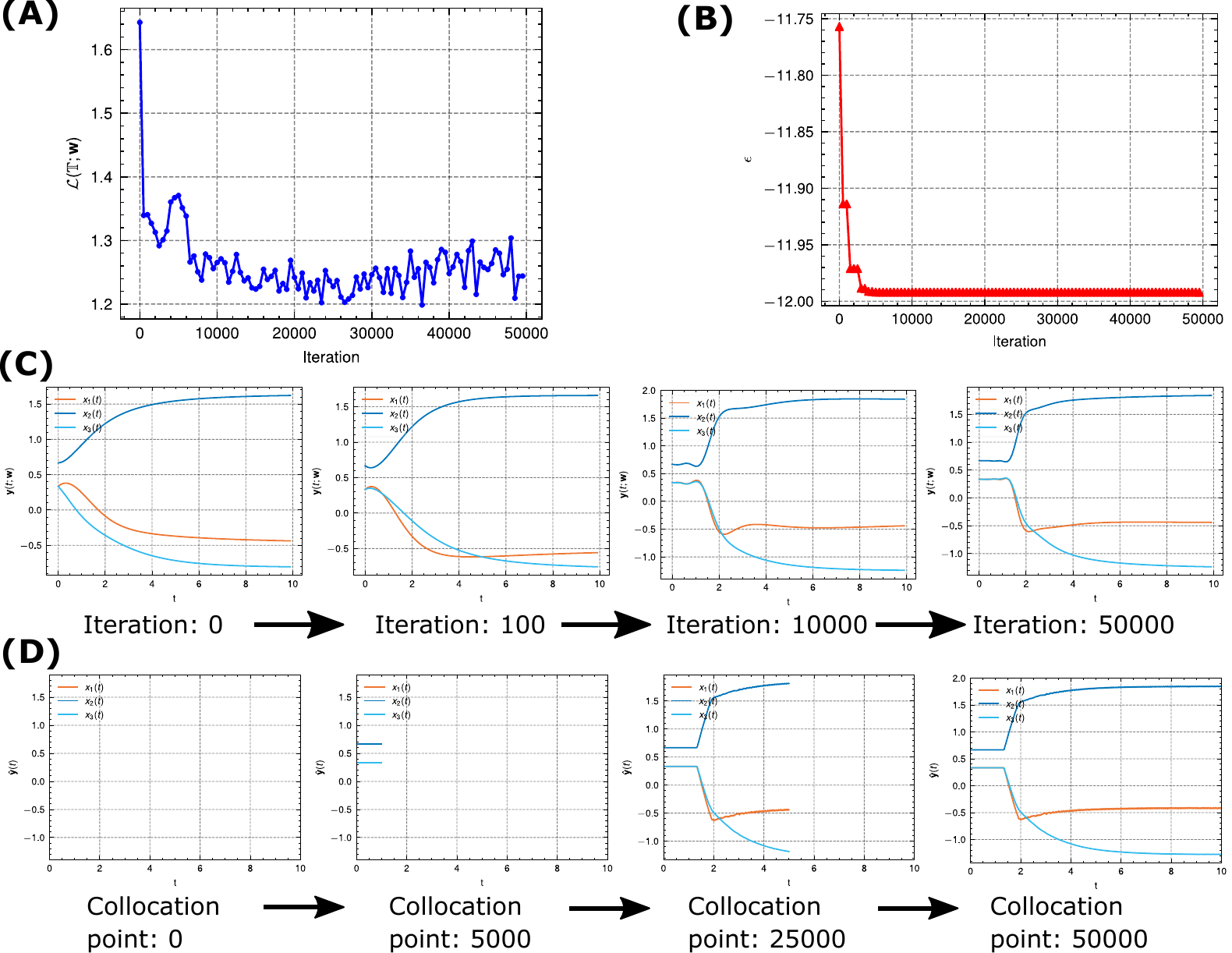}
    \caption{\textbf{Example 6 pseudoconvex nonsmooth standard CNLP} 
    \textbf{(A) The loss versus the number of iterations.} 
    \textbf{(B) The epsilon value versus the number of iterations.} The epsilon metric is defined in \eqref{eq: epsilon 2}
    \textbf{(C) The solving process of the OINN model}
    \textbf{(D) The solving process of the  numerical integration method}}
    \label{fig: Example 6}
\end{figure}

\begin{table}[!ht]
\begin{tabular}{@{}l|ll|ll@{}}
\toprule
\multirow{2}{*}{Index}     & \multicolumn{2}{l|}{OINN}                       & \multicolumn{2}{l}{Numerical integration method}      \\ \cmidrule(l){2-5} 
                           & Iteration & Solution                            & Collocation point & Solution \\ \midrule
\multirow{6}{*}{Example 6} & 0         & [-0.44,  1.62, -0.81]               & 0                 & [0.33, 0.67, 0.33]         \\
                           & 10        & [-0.44,  1.62, -0.81]               & 10                & [0.33, 0.67, 0.33]        \\
                           & 100       & [-0.56,  1.66, -0.76]               & 100               &  [0.33, 0.67, 0.33]        \\
                           & 1000      & [-0.52,  1.70, -0.89]               & 1000              &  [0.33, 0.67, 0.33]        \\
                           & 10000     & [-0.44,  1.84, -1.24]               & 10000             &  [-0.63, 1.56, -0.49]        \\
                           & 50000     & [-0.44,  1.84, -1.24]               & 50000             &  [-0.41, 1.85, -1.28]       \\ \bottomrule
\end{tabular}
\centering
\caption{\textbf{Example 6, Approximate solutions to the NPE during solving}}
\label{tab: Example 6}
\end{table}

\DW{Figure  \ref{fig: Example 6} shows the training of the OINN model, where the loss decreased from $1.64$ to $1.20$, and the epsilon value decreased from $-11.75$ to $-11.99$. Figure \ref{fig: Example 6} (C) and (D) show the progressions of the approximate state solutions to the IVP \eqref{eq: IVP-Example6}.}

\DW{Table \ref{tab: Example 6} displays the progressions of the approximate solutions to the CNLP \eqref{eq: Example 6}.
The OINN model gives the final solution of $[-0.44, 1.84, -1.24]$, whereas the numerical integration method give the final solution of $[-0.41, 1.85, -1.28]$.}

\subsection{Hyperparameters study}\label{section: 5.2}

\DW{In this subsection, we discuss the setting of the two critical hyperparameters in OINN, i.e., initial point and time range. The experiments are based on Example 3 of Section \ref{section: 5.1.3} \AL{for illustrative purposes}.}

\begin{figure}[!ht]
    \centering
    \includegraphics{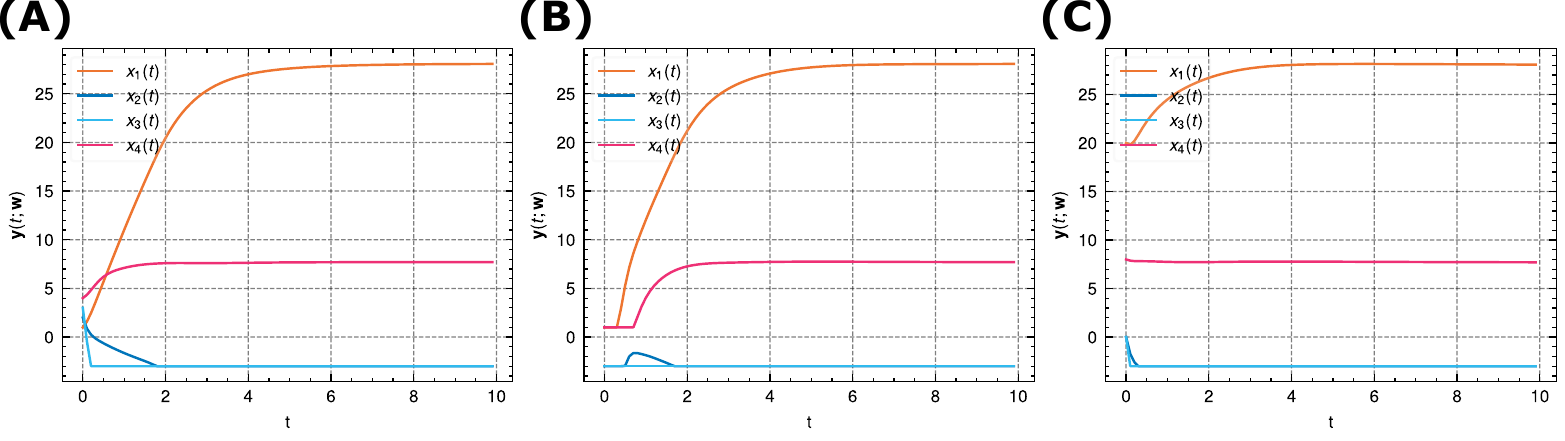}
    \caption{\textbf{The OINN models with different initial points} 
    \textbf{(A) Initial point $\mathbf{y}_0 = [1, 2, 3, 4]$} 
    \textbf{(B) Initial point $\mathbf{y}_0 = [-10, -15, -10, -14]$}
    \textbf{(C) Initial point $\mathbf{y}_0 = [20, 0, 0, 8]$} }
    \label{fig: Hyper initial point}
\end{figure}

\begin{table}[!ht]
\begin{tabular}{lllll}
\hline
Time range                             &Initial point                             & Iteration & Solution                            & epsilon $\downarrow$ \\ \hline
\multirow{15}{*}{$t\in[0,10]$}         &\multirow{5}{*}{$\mathbf{y}_0 = [1, 2, 3, 4]$}             & 0         & [1.46, 2.20,  3.49, 4.21]           &  6.49    \\
                                       &                                          & 10        & [2.21, 1.46, 2.66, 4.80 ]           &  5.66   \\
                                       &                                          & 100       & [2.53, 1.17, 2.32, 4.97]            & 5.32    \\
                                       &                                          & 1000      & [25.65, -3.00,   -3.00,    8.00]    & 2.42    \\
                                       &                                          & 10000     & [28.07, -3.00,   -3.00,    7.71]    &  0.00    \\ \cline{2-5}
                                       &\multirow{5}{*}{$\mathbf{y}_0 = [-10, -15, -10, -14]$}     & 0         & [ 1.00, -3.00, -3.00,  1.00]        & 28.11        \\
                                       &                                          & 10        & [ 1.00, -3.00, -3.00,  1.00]        & 28.11        \\
                                       &                                          & 100       & [ 1.00, -3.00, -3.00,  1.00]        & 28.11        \\
                                       &                                          & 1000      & [15.93, -1.35, -3.00, 9.69]         & 4.04        \\
                                       &                                          & 10000     & [28.07, -3.00, -3.00, 7.71]         & 0.00         \\ \cline{2-5}
                                       &\multirow{5}{*}{$\mathbf{y}_0 = [20, 0, 0, 8]$}            & 0         & [18.92, -0.92,  0.24,  7.53]        & 3.24        \\
                                       &                                          & 10        & [18.52, -1.58, -0.46,  6.88]        & 2.54        \\
                                       &                                          & 100       & [18.52, -1.58, -0.46,  6.88]        & 2.54        \\
                                       &                                          & 1000      & [28.30, -3.00,  -3.00,   7.6]       & 0.23        \\
                                       &                                          & 10000     & [28.07, -3.00,  -3.00,  7.71]       & 0.00        \\ \hline
\end{tabular}
\centering
\caption{\textbf{OINN solutions with different initial points}}
\label{tab: Hyper initial point}
\end{table}

\DW{\textbf{Initial point} \ Figure \ref{fig: Hyper initial point} and Table \ref{tab: Hyper initial point} show the convergence behavior of three different initial points. 
Thanks to the global convergence property of the ODE system, any initial point can converge to the optimal solution, provided there are large enough training iterations. 
The convergence occurs faster and requires fewer training iterations if the initial point is closer to the optimal solution.
For example, the initial point $[20, 0, 0, 8]$ is the closest to the optimal solution $[28.07, -3.00, -3.00, 7.71]$ among the three, so it reaches the lowest epsilon value of $0.23$ at the 1000th iteration, while the epsilon values of the other two are $2.42$ and $4.04$.}

\begin{figure}[!ht]
    \centering
    \includegraphics{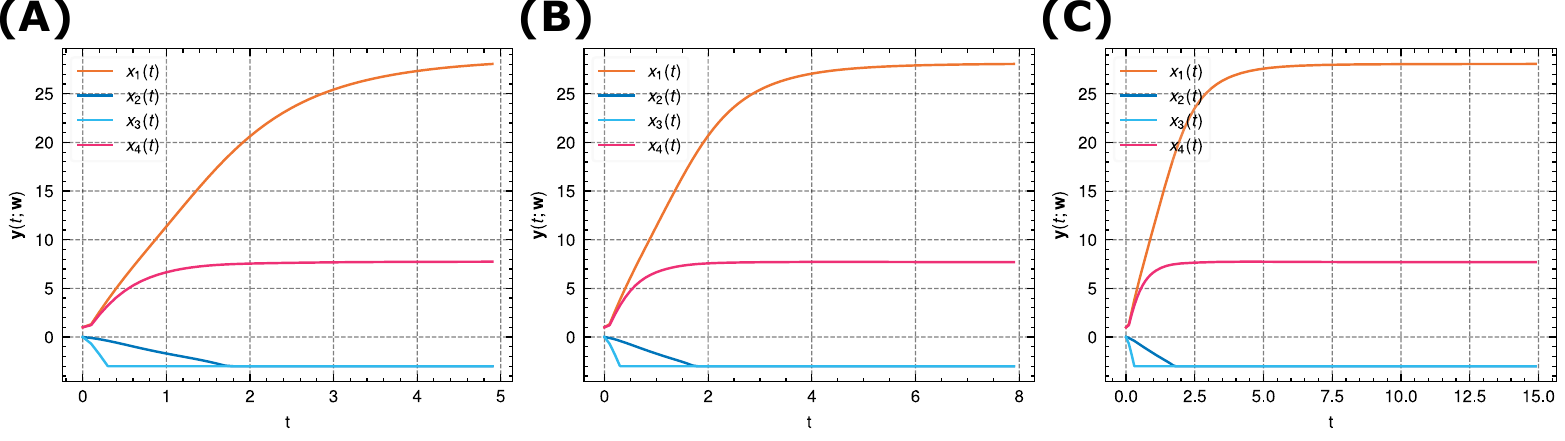}
    \caption{\textbf{The OINN models with different time ranges} 
    \textbf{(A) Time range $t\in[0, 5]$} 
    \textbf{(B) Time range $t\in[0, 8]$} 
    \textbf{(C) Time range $t\in[0, 15]$} }
    \label{fig: Hyper time range}
\end{figure}

\begin{table}[!ht]
\begin{tabular}{lllll}
\hline
Initial point                                       & Time range               & Iteration & Solution & epsilon $\downarrow$ \\ \hline
\multirow{15}{*}{$\mathbf{y}_0 = [0, 0, 0, 0]$}    & \multirow{5}{*}{$t\in[0, 5]$}  & 0         & [ 1.00,   -0.71, -0.42,  1.00]         & 16.66        \\
                                                    &                          & 10        & [ 1.00,   -0.59, -1.24,  1.36]         & 16.08       \\
                                                    &                          & 100       & [ 3.77, -0.25, -3.00,    5.30]         & 10.92        \\
                                                    &                          & 1000      & [25.82, -3.00,   -3.00,    7.92]         & 2.26        \\
                                                    &                          & 10000     & [28.11, -3.00,   -3.00,    7.74]         & 0.07        \\ \cline{2-5}
                                                    & \multirow{5}{*}{$t\in[0, 8]$}  & 0         & [1.00,   1.47, 2.12, 1.48]         &  8.08       \\
                                                    &                          & 10        & [1.00,   1.10,  1.74, 1.83]         &  7.64       \\
                                                    &                          & 100       & [1.00,   1.10,  1.74, 1.83]         &  7.64       \\
                                                    &                          & 1000      & [26.42, -3.00,   -3.00,    8.13]        &  1.65       \\
                                                    &                          & 10000     & [28.07, -3.00,   -3.00,    7.71]         &  0.00       \\ \cline{2-5}
                                                    & \multirow{5}{*}{$t\in[0, 15]$} & 0         & [ 1.00,    0.02, -0.13,  1.00  ]         &  12.99       \\
                                                    &                          & 10        & [ 1.00,    0.02, -0.13,  1.00  ]         &  12.99       \\
                                                    &                          & 100       & [ 3.07,  0.15, -3.00,    3.50 ]         &  9.69       \\
                                                    &                          & 1000      & [25.98, -3.00,   -3.00,    7.76]        & 2.09        \\
                                                    &                          & 10000     & [28.07, -3.00,   -3.00,    7.70]         &  0.00       \\ \hline
\end{tabular}
\centering
\caption{\textbf{OINN solutions with different time ranges}}
\label{tab: Hyper time range}
\end{table}

\DW{\textbf{Time range} \ 
Figure \ref{fig: Hyper time range} and Table \ref{tab: Hyper time range} show the OINN model of the same initial point with different time ranges.
The first time range $[0, 5]$ still has an epsilon value of $0.07$ after 10,000 iterations. This is because the OINN model has reached the upper limit of accuracy.
The time ranges $[0, 8]$ and $[0, 15]$ both eventually reach an epsilon value of $0.00$. 
The time range $[0, 8]$ converges faster than the time range $[0, 15]$ because the former has a smaller span and is easier to train \AL{than the latter}.}

\subsection{Discussions}\label{section: 5.3}

\DW{The OINN method and numerical integration methods solve CNLP in different ways and are based on different software implementations.
Because of \AL{that}, it is difficult to determine which method is superior to the other.
In this subsection, we first highlight some of OINN's computational advantages and then discuss its limitations.}

\begin{table}[!ht]
\begin{tabular}{l|lc|lc}
\hline
\multirow{2}{*}{Index}     & \multicolumn{2}{l|}{OINN} & \multicolumn{2}{l}{Numerical integration method} \\ \cline{2-5} 
                           & Iteration    & Epsilon: NPE error  $\downarrow$  & Collocation point      & Epsilon: NPE error $\downarrow$    \\ \hline
\multirow{6}{*}{Example 3} & 0            & 10.925     & 0                      & 13.111      \\
                           & 10           & 10.925     & 10                     & 13.123      \\
                           & 100          & 10.925     & 100                    & 13.223      \\
                           & 1000         & \textbf{1.550}      & 1000                   & 12.093      \\
                           & 10000        & \textbf{0.001}      & 10000                  & 7.307       \\
                           & 50000        & \textbf{0.001}      & 50000                  & 0.006       \\ \hline
\multirow{6}{*}{Example 4} & 0            & 3.117      & 0                      & 2.437       \\
                           & 10           & 0.738      & 10                     & 2.350       \\
                           & 100          & \textbf{0.007}      & 100                    & 1.728       \\
                           & 1000         & \textbf{0.000}      & 1000                   & 0.160       \\
                           & 10000        & \textbf{0.000}      & 10000                  & 0.001       \\
                           & 50000        & \textbf{0.000}      & 50000                  & 0.001       \\ \hline
\end{tabular}
\centering
\caption{\textbf{Comparison of solutions accuracy} The epsilon metric is defined in \eqref{eq: epsilon 1}.}
\label{tab: pro 1.1}
\end{table}

\begin{table}[!ht]
\begin{tabular}{l|lc|lc}
\hline
\multirow{2}{*}{Index}     & \multicolumn{2}{l|}{OINN}   & \multicolumn{2}{l}{Numerical integration method} \\ \cline{2-5} 
                           & Iteration & Epsilon: Objective value $\downarrow$ & Collocation point  & Epsilon:  Objective value $\downarrow$ \\ \hline
\multirow{6}{*}{Example 5} & 0         & 62.315          & 0                  & 43.838           \\
                           & 10        & 39.633          & 10                 & 43.757          \\
                           & 100       & 39.020          & 100                & 43.130          \\
                           & 1000      & 39.020          & 1000               & 40.120          \\
                           & 10000     & 39.020          & 10000              & 39.029           \\
                           & 50000     & 39.020          & 50000              & 39.020           \\ \hline
\multirow{6}{*}{Example 6} & 0         & -11.757         & 0                  & -7.871          \\
                           & 10        & -11.757         & 10                 & -7.871          \\
                           & 100       & -11.861         & 100                & -7.871          \\
                           & 1000      & -11.914         & 1000               & -7.871          \\
                           & 10000     & \textbf{-11.992}         & 10000              & inf             \\
                           & 50000     & \textbf{-11.992}         & 50000              & -11.985         \\ \hline
\end{tabular}
\centering
\caption{\textbf{Comparison of objective values} The epsilon metric is defined in \eqref{eq: epsilon 2}. inf means that the solution is not feasible.}
\label{tab: pro 1.2}
\end{table}

\DW{Tables \ref{tab: pro 1.1} and \ref{tab: pro 1.2} display the epsilon values \AL{ while Examples 3-6 were being resolved}, where the OINN model converges with training iterations and the numerical integration method converges as the collocation point \AL{progresses}.
In Examples 3 and 4, the OINN model reaches a lower final epsilon error than the numerical integration method, indicating that OINN found a better solution to satisfy the CNLPs.
In Example 6, the OINN model reaches a lower objective value of $-11.992$, while the numerical integration method only manages to reach the objective value of $-11.985$.}

\DW{OINN can give an approximate solution to the CNLP at any round of iterations, while the numerical integration method can only give the solution at the end of the program. 
Thanks to that, OINN can provide more accurate approximate solutions in the early stage of \AL{the solving process}.
For instance, in Example 3, the OINN's epsilon error has decreased to $1.55$ by the 1000th training iteration, compared to $12.093$ of the numerical integration method.}

\begin{table}[!ht]
\begin{adjustbox}{width=1\textwidth}
\begin{tabular}{l|ll|lllllll}
\hline
\multirow{2}{*}{Index}     & \multicolumn{2}{l|}{OINN} & \multicolumn{7}{l}{Numerical integration methods}                                                                                                                                                                                                                                                                                                                                                                                                   \\ \cline{2-10} 
                           & Iteration    & CPU time   & \begin{tabular}[c]{@{}l@{}}Collocation\\ point\end{tabular} & \begin{tabular}[c]{@{}l@{}}RK45\\ CPU\\ time\end{tabular} & \begin{tabular}[c]{@{}l@{}}RK23\\ CPU\\ time\end{tabular} & \begin{tabular}[c]{@{}l@{}}DOP853\\ CPU\\ time\end{tabular} & \begin{tabular}[c]{@{}l@{}}Radau\\ CPU\\ time\end{tabular} & \begin{tabular}[c]{@{}l@{}}BDF\\ CPU\\ time\end{tabular} & \begin{tabular}[c]{@{}l@{}}LSODA\\ CPU\\ time\end{tabular} \\ \hline
\multirow{5}{*}{Example 6} & 10           & 202 ms     & 10                                                          & 1350 ms                                                     & 860 ms                                                     & 3470 ms                                                      & 1000 ms                                                        & Fail                                                     & 157 ms                                                     \\
                           & 100          & 893 ms     & 100                                                         & 1740 ms                                                     & 1090 ms                                                     & 4620 ms                                                      & 1330 ms                                                     & Fail                                                     & 154 ms                                                     \\
                           & 1000         & 8.47 s     & 1000                                                        & 2.14s                                                     & 1.32s                                                     & 5.68 s                                                      & 1.47 s                                                     & Fail                                                     & 188 ms                                                     \\
                           & 10000        & 1min 20s       & 10000                                                       & 1min 25s                                                  & 5min 5s                                                   & 34min 29s                                                   & Fail                                                       & Fail                                                     & 4h 4min 35s                                                \\
                           & 50000        & \textbf{7min 55s}      & 50000                                                       & 14min 29s                                                 & 25min 14s                                                 & 1h 43min 28s                                                & Fail                                                       & Fail                                                     & Fail                                                       \\ \hline
\end{tabular}
\end{adjustbox}
\centering
\caption{\textbf{CPU times of the OINN method and numerical integration methods.} RK45, RK23, DOP853, Radau, BDF and LSODA are six different numerical integration methods. ms, s, min, and h refer to milliseconds, seconds, minutes, and hours respectively.}
\label{tab: pro 2}
\end{table}

\DW{Table \ref{tab: pro 2} shows the CPU times for solving Example 6, which has a stiff ODE system and is challenging to solve numerically.
We compare OINN  with six different numerical integration methods.
OINN outperforms all these six methods in terms of computational \AL{CPU time, i.e.,} OINN takes 7min 55s while RK45 \AL{takes at best} 14min 29s. The three methods, Radau, BDF, and LSODA, fail to solve this problem.}

\DW{We must emphasize that the proposed OINN should not be seen as a substitute for the conventional numerical integration methods.
Such methods have been developed for many years and are well known to meet the requirements of reliability.
OINN research is still in its early stages, making it difficult to go beyond traditional methods for many practical problems.
Our contribution is to open up a fresh perspective and a new line of research \AL{to solve CNLP}.}

\section{Conclusion}\label{section: 6}
% Summary of the paper
We propose a deep learning approach to address CNLPs, namely OINN.
We give a complete description of solving CNLPs by OINN, including neurodynamic optimization approaches, the OINN framework, and the training algorithm. By doing so, we connect this longstanding problem with deep learning and machine learning communities. With the rapid development of deep learning research, both methodologically and experimentally, we believe that this work will lead to ongoing contributions that can benefit a wide range of practitioners in optimization.

% Future directions
\DW{There are many possible future directions for this work; we give some examples here. 
1) How to design a better neural network structure and activation function for the OINN model? 
2) What is the difference \AL{terms of computational effort and quality of the solution} when using various epsilon metrics for the same CNLP?
3) How \AL{to design appropriate epsilon} metrics for other CNLPs?
%4) How to apply OINN to more complex optimization problems? 
We can gradually make our proposed approach more robust by providing answers to these questions.}
\bibliography{OINN}
\end{document}